\documentclass{amsart}
\usepackage{amssymb}
\usepackage{amsmath}
\usepackage{amsfonts}

\newtheorem{theorem}{Theorem}
\theoremstyle{plain}

\newtheorem{corollary}{Corollary}

\newtheorem{definition}{Definition}

\newtheorem{lemma}{Lemma}

\newtheorem{proposition}{Proposition}
\newtheorem{remark}{Remark}

\numberwithin{equation}{section}
\input{tcilatex}

\begin{document}
\title[]{Notes concerning K\"{a}hler and anti-K\"{a}hler structures on
quasi-statistical manifolds }
\author{Aydin GEZER}
\address{Ataturk University, Faculty of Science, Department OF Mathematics,
25240, Erzurum-Turkey.}
\email{aydingzr@gmail.com}
\thanks{}
\author{Bu\c{s}ra AKTA\c{S}}
\address{Kirikkale University, Faculty of Science and Arts, Department OF
Mathematics, 71450, Kirikkale-Turkey}
\email{baktas6638@gmail.com}
\urladdr{}
\thanks{}
\author{Olgun DURMAZ}
\address{Kyrgyz-Turkish Manas University, Faculty of Science, Department of
Mathematics, Bishkek, Kyrgyzstan}
\email{durmazolgun@gmail.com}
\urladdr{}
\date{May 26, 2023}
\thanks{This paper is in final form and no version of it will be submitted
for publication elsewhere.}

\begin{abstract}
Let $(\acute{N},g,\nabla )$\ be a $2n$-dimensional quasi-statistical
manifold that admits a pseudo-Riemannian metric $g$ (or $h)$ and a linear
connection $\nabla $ with torsion. This paper aims to study an almost
Hermitian structure $(g,\tciLaplace )$ and an almost anti-Hermitian
structure $(h,\tciLaplace )$ on a quasi-statistical manifold that admit an
almost complex structure $\tciLaplace $. Firstly, under certain conditions,
we present the integrability of the almost complex structure $\tciLaplace $.
We show that when $d^{\nabla }\tciLaplace =0$ and the condition of
torsion-compatibility are satisfied, $(\acute{N},g,\nabla ,$ $\tciLaplace )$
turns into a K\"{a}hler manifold. Secondly, we give necessary and sufficient
conditions under which $(\acute{N},h,\nabla ,\tciLaplace )$ is an anti-K\"{a}%
hler manifold, where $h$ is an anti-Hermitian metric. Moreover, we search
the necessary conditions for $(\acute{N},h,\nabla ,\tciLaplace )$ to be a
quasi-K\"{a}hler-Norden manifold.


\end{abstract}

\maketitle

\section{Introduction}

Nowadays, the research of spaces that consists of probability measures is
getting more attention. Information geometry that is a famous theory in
geometry is a tool to search such spaces. Information geometry as a theory
which combines differential geometry and statistics has a great importance
in science. For example; image processing, physics, computer science and
machine learning are some of its applications (\cite%
{Amari,Belkin,Caticha,Sun}). Two geometric quantities, called dual
connections describing the derivation with respect to vector fields play an
important role in characterizing statistical manifolds. The investigation of
dual elements and the relationships between them constitute the basic
direction of advancement in the research of statistical manifolds \cite%
{Calin}. Recently, the notion of statistical manifold has attracted the
attention of many mathematicians (\cite{Balan,Fei,Furuhata,Peyghan,Teofilova}%
).

The notion of statistical manifold admitting torsion (SMAT) (or quasi
statistical manifolds) in order to define geometric structures on quantum
state spaces was firstly introduced by Kurose \cite{Kurose}. This manifold
naturally appears in a quantum statistical model and can be regarded as the
quantum version of statistical manifolds. A statistical manifold admitting
torsion (SMAT) is a pseudo Riemannian manifold with a pair of dual
connections, where only one of them must be torsion free but the other is
not necessarily so. The expression of SMAT was originally represented to
study such a geometrical structure from a mathematical point of view (\cite%
{Amari1,Kurose}).

Let $\acute{N}$ be a $2n$-dimensional differentiable manifold and $g$ be a
pseudo Riemannian metric on $\acute{N}$. An almost complex structure on $%
\acute{N}$ is a tensor field $\tciLaplace $ of type $\left( 1,1\right) $
such that $\tciLaplace ^{2}=-id$. An almost complex manifold is such a
manifold with a fixed almost complex structure. Note that almost complex
structures exist only when $\acute{N}$ is of even dimension. Ensuring
compatibility of $\tciLaplace $ with $g$, $g\left( \tciLaplace \xi _{1},\xi
_{2}\right) +g\left( \xi _{1},\tciLaplace \xi _{2}\right) =0$, for any
vector fields $\xi _{1}$ and $\xi _{2}$ on $\acute{N}$, leads to an almost
Hermitian manifold $\left( \acute{N},g,\tciLaplace \right) $. The compatible
metric is also called a Hermit metric. If $\tciLaplace $ is integrable, the
manifold $\left( \acute{N},g,\tciLaplace \right) $ becomes a Hermitian
manifold. Moreover, the fundamental two form $\omega $ can be described $%
\omega \left( \xi _{1},\xi _{2}\right) =g\left( \tciLaplace \xi _{1},\xi
_{2}\right) $ and performs to satisfy $\omega \left( \tciLaplace \xi
_{1},\xi _{2}\right) +\omega \left( \xi _{1},\tciLaplace \xi _{2}\right) =0$%
. An almost K\"{a}hler manifold is an almost Hermitian manifold whose
fundamental $2$-form $\omega $ is closed. In other words, an almost K\"{a}%
hler manifold is a symplectic manifold equipped with a compatible metric.
With integrability of $\tciLaplace $, the almost K\"{a}hler manifold $\left( 
\acute{N},g,\tciLaplace \right) $ rises to a K\"{a}hler manifold. Also, it
is well known that the almost Hermitian manifold $\left( \acute{N}%
,g,\tciLaplace \right) $ is K\"{a}hler manifold if and only if the almost
complex structure $\tciLaplace $ is covariantly constant with respect to the
Levi-Civita connection $\nabla ^{g}$, that is, $\nabla ^{g}\tciLaplace =0$.
\ Fei and Zhang presented an alternative classification for K\"{a}hler
manifolds by taking any torsion-free linear connection $\nabla $ instead of
the Levi-Civita connection $\nabla ^{g}$ \cite{Fei}. Firstly, they showed
that Codazzi coupling of a torsion-free linear connection $\nabla $ with $%
\tciLaplace $ implies the integrability of $\tciLaplace $. Furthermore, they
proved that a torsion-free linear connection $\nabla $ is Codazzi-coupled
with both $g$ and $\tciLaplace $, then the triple $\left( \acute{N}%
,g,\tciLaplace \right) $ turns into a K\"{a}hler manifold. Such a K\"{a}hler
manifold is called Codazzi K\"{a}hler manifold \cite{Fei}.

An anti-K\"{a}hler (Norden-K\"{a}hler) manifold means a manifold $\left( 
\acute{N},h,\tciLaplace \right) $ which consists of a differentiable
manifold of dimension $2n$, an almost complex structure $\tciLaplace $ and
an anti-Hermitian metric $h$ such that $\nabla ^{h}\tciLaplace =0$, where $%
\nabla ^{h}$ is the Levi-Civita connection of $h$ \cite{Ganchev,Manev}. The
metric $h$ is called an anti-Hermitian (Norden) metric if it satisfies $%
h\left( \tciLaplace \xi _{1},\xi _{2}\right) -h\left( \xi _{1},\tciLaplace
\xi _{2}\right) =0$ for all vector fields $\xi _{1}$ and $\xi _{2}$ on $%
\acute{N}$. Then the metric $h$ has necessarily a neutral signature $(n,n)$.
By $\hslash (\xi _{1},\xi _{2})=h\left( \tciLaplace \xi _{1},\xi _{2}\right) 
$, the twin metric $\hslash $ can be defined and it is symmetric and
satisfies $\hslash \left( \tciLaplace \xi _{1},\xi _{2}\right) -\hslash
\left( \xi _{1},\tciLaplace \xi _{2}\right) =0$ for any vector fields $\xi
_{1},\xi _{2}$ on $\acute{N}$, that is, it is another anti-Hermitian
(Norden) metric. Since there exists a pair of anti-Hermitian (or Norden)
metrics on anti-Hermitian (or Norden) manifolds, one can take into
consideration dual (conjugate) connections according to each of these metric
tensors and their relations to dual connections related to the almost
complex structure. Hence, the investigation of statistical structures on
these manifolds is of great importance. Hermitian manifolds as well as
Norden manifolds have been analyzed from various points of view. Here, we
refer to (\cite{Fei,Gezer,Grigorian}). In \cite{Iscan}, the authors
presented a new approach for expanding almost anti-Hermitian manifolds to
anti-K\"{a}hler manifolds. They showed that the anti-K\"{a}hler condition is
equivalent to $%
\mathbb{C}
$ analyticity of the anti-Hermitian metric $h$, that is, $\Phi _{\tciLaplace
}h=0$, where $\Phi _{\tciLaplace }$ is the Tachibana operator. Moreover,
using Codazzi coupling of $\left( \nabla ,\tciLaplace \right) $, Gezer and
Cakicioglu gave an alternative classification for anti-K\"{a}hler manifolds
with respect to a torsion-free linear connection $\nabla $ \cite{Gezer}.
After then, considering the presence of Tachibana operator and Codazzi
coupling of $\left( \nabla ,g\right) $ with a torsion-free linear connection 
$\nabla $, Durmaz and Gezer have showed that locally metallic
pseudo-Riemannian manifolds can be classified \cite{Durmaz}.

Now, it is natural to ask the following question: \textit{Can the K\"{a}hler
and anti-K\"{a}hler manifolds be classified by taking the any linear
connection }$\nabla $\textit{\ with torsion tensor }$T^{\nabla }$\textit{\
instead of the Levi-Civita connection }$\nabla ^{g}$\textit{\ }$\left( \text{%
resp. }\nabla ^{h}\right) $\textit{\ of }$g$\textit{\ }$\left( \text{resp. }%
h\right) $? We can also ask this question as the following: \textit{Does the
torsion tensor of the linear connection }$\nabla $\textit{\ have to be flat
in order to classify these manifolds?} This paper aims to find the answers
to these questions. The present paper is organized as follows: In section 2,
for any linear connection $\nabla $ with torsion tensor $T^{\nabla }$, we
investigate the integrability of the almost complex structure $\tciLaplace $
in Lemma \ref{lem1} and Proposition \ref{pro2}. By using the definitions of $%
d^{\nabla }\tciLaplace $ and the Vishnevskii operator $\Psi _{\tciLaplace }$%
, we get different results between the integrability of $\tciLaplace $ and $%
d^{\nabla }\tciLaplace $. In section 3, considering the definition of
quasi-statistical structure $\left( \nabla ,g\right) $ with a
pseudo-Riemannian metric $g$, we obtain interesting relationships between
the quasi structures of conjugate connections $\nabla ^{\ast },\nabla
^{\dagger }$ and $\nabla ^{\tciLaplace }$, and $d^{\nabla \left( \nabla
^{\ast },\nabla ^{\dagger },\nabla ^{\tciLaplace }\right) }$-closed of $%
\tciLaplace $ $\left( \text{see Proposition \ref{pro5}}\right) $. After
then, we show that the triple $\left( \acute{N},g,\tciLaplace \right) $ is a
K\"{a}hler manifold if $d^{\nabla }\tciLaplace =0$ and $T^{\nabla }\left(
\tciLaplace \xi _{1},\xi _{2}\right) =-T^{\nabla }\left( \xi
_{1},\tciLaplace \xi _{2}\right) $ on a quasi-statistical manifold $\left( 
\acute{N},g,\nabla \right) $, where $\nabla $ is any linear connection with
torsion tensor $T^{\nabla }$ (see Theorem \ref{teo2}). In the last section,
by taking an anti-Hermitian metric $h$ instead of a Hermitian metric $g$, we
reinvestigate the properties of quasi-statistical structures. Taking any
linear connection $\nabla $ with torsion tensor $T^{\nabla }$ instead of $%
\nabla ^{h}$ which is the Levi-Civita connection of $h$, for anti-K\"{a}hler
and quasi-K\"{a}hler-Norden manifolds, we showed that the new
classifications can be constituted (see Theorem \ref{teo5} and Theorem \ref%
{cor8}).

\section{$\tciLaplace $-Conjugate of $\bigtriangledown $ and $%
d^{\bigtriangledown }$-closed of $\tciLaplace $}

In this section, we study a linear connection on a differentiable manifold $%
\acute{N}$ with a $(1,1)$-tensor field $\tciLaplace $, which satisfies the
condition $d^{\bigtriangledown }\tciLaplace =0$. The structure $\tciLaplace $
is called an almost complex structure if $\tciLaplace ^{2}=-id$ or an almost
paracomplex structure if $\tciLaplace ^{2}=id$. In this paper, we will deal
with the almost complex structure $\tciLaplace $. Note that the similar
results may be obtained if you use the almost para complex structure instead
of the almost complex structure.\newline

Starting from a linear connection $\bigtriangledown $ on $\acute{N}$, we can
apply an $\tciLaplace $-conjugate transformation to achieve a new connection 
$\bigtriangledown ^{\tciLaplace }:=\tciLaplace ^{-1}\bigtriangledown
\tciLaplace $ or $\bigtriangledown _{\xi _{1}}^{\tciLaplace }\xi
_{2}=\tciLaplace ^{-1}\left( \bigtriangledown _{\xi _{1}}\tciLaplace \xi
_{2}\right) $ for any vector fields $\xi _{1}$ and $\xi _{2}$, where $%
\tciLaplace ^{-1}$ identifies the inverse isomorphism of $\tciLaplace $. It
can be confirmed that indeed $\bigtriangledown ^{\tciLaplace }$ is a linear
connection.

\begin{definition}
A linear connection $\nabla $ and a $(1,1)$-tensor field $\tciLaplace $ are
called Codazzi-coupled if the following equality exists%
\begin{equation*}
\left( \nabla _{\xi _{1}}\tciLaplace \right) \xi _{2}=\left( \nabla _{\xi
_{2}}\tciLaplace \right) \xi _{1},
\end{equation*}%
where $\left( \nabla _{\xi _{1}}\tciLaplace \right) \xi _{2}=\nabla _{\xi
_{1}}\tciLaplace \xi _{2}-\tciLaplace \nabla _{\xi _{1}}\xi _{2}$.
\end{definition}

As a linear connection, $\bigtriangledown $ creates a map $\bigtriangledown
:\Omega ^{0}\left( T\acute{N}\right) \rightarrow \Omega ^{1}\left( T\acute{N}%
\right) $, where $\Omega ^{i}\left( T\acute{N}\right) $ is the space of
smooth $i$-forms with value in $T\acute{N}$. Regarding $\tciLaplace $ as an
element of $\Omega ^{1}\left( T\acute{N}\right) ,$ it is easy to see 
\begin{equation*}
\left( d^{\bigtriangledown }\tciLaplace \right) \left( \xi _{1},\xi
_{2}\right) =\left( \bigtriangledown _{\xi _{1}}\tciLaplace \right) \xi
_{2}-\left( \bigtriangledown _{\xi _{2}}\tciLaplace \right) \xi
_{1}+\tciLaplace T^{\bigtriangledown }\left( \xi _{1},\xi _{2}\right) ,
\end{equation*}%
where the torsion tensor is given by $T^{\bigtriangledown }\left( \xi
_{1},\xi _{2}\right) =\bigtriangledown _{\xi _{1}}\xi _{2}-\bigtriangledown
_{\xi _{2}}\xi _{1}-\left[ \xi _{1},\xi _{2}\right] $. Hence, $\tciLaplace $
is called $d^{\bigtriangledown }$-closed if $d^{\bigtriangledown
}\tciLaplace =0$.

Any $(1,1)$-tensor field $\tciLaplace $ is called a quadratic operator if
there exists $\alpha \neq \beta \in 
\mathbb{C}
$ such that $\alpha +\beta $ and $\alpha \beta $ are real numbers and $%
\tciLaplace ^{2}-\left( \alpha +\beta \right) \tciLaplace +\alpha \beta
.id=0.$ Note that $\tciLaplace $ is an isomorphism, so $\alpha \beta \neq 0$.

The Nijenhuis tensor $N_{\tciLaplace }$ associated with $\tciLaplace $ is
described as%
\begin{equation*}
N_{\tciLaplace }\left( \xi _{1},\xi _{2}\right) =-\tciLaplace ^{2}\left[ \xi
_{1},\xi _{2}\right] +\tciLaplace \left[ \xi _{1},\tciLaplace \xi _{2}\right]
+\tciLaplace \left[ \tciLaplace \xi _{1},\xi _{2}\right] -\left[ \tciLaplace
\xi _{1},\tciLaplace \xi _{2}\right] \text{.}
\end{equation*}%
When $N_{\tciLaplace }=0$, $\tciLaplace $ is said to be integrable.

\begin{proposition}
\label{GAD1}Let $\bigtriangledown $ be a linear connection and let $%
\tciLaplace $ be a $(1,1)$-tensor field on $\acute{N}$. Then, there exist
the following expressions

$\left( i\right) $ $d^{\bigtriangledown }\tciLaplace =0\Leftrightarrow
T^{\bigtriangledown ^{\tciLaplace }}=0$;

$\left( ii\right) $ $d^{\bigtriangledown ^{\tciLaplace }}\tciLaplace
=0\Leftrightarrow T^{\bigtriangledown }=0$;

$\left( iii\right) $ $d^{\bigtriangledown }\tciLaplace =d^{\bigtriangledown
^{\tciLaplace }}\tciLaplace \Leftrightarrow \left( \bigtriangledown
,\tciLaplace \right) $ is Codazzi-coupled.
\end{proposition}

\begin{proof}
$\left( i\right) $ For any vector fields $\xi _{1}$ and $\xi _{2}$, we have%
\begin{eqnarray*}
&&\left( d^{\bigtriangledown }\tciLaplace \right) \left( \xi _{1},\xi
_{2}\right) \\
&=&\left( \bigtriangledown _{\xi _{1}}\tciLaplace \right) \xi _{2}-\left(
\bigtriangledown _{\xi _{2}}\tciLaplace \right) \xi _{1}+\tciLaplace
T^{\bigtriangledown }\left( \xi _{1},\xi _{2}\right) \\
&=&\bigtriangledown _{\xi _{1}}\tciLaplace \xi _{2}-\bigtriangledown _{\xi
_{2}}\tciLaplace \xi _{1}-\tciLaplace \left[ \xi _{1},\xi _{2}\right] \\
&=&\tciLaplace \left( \tciLaplace ^{-1}\bigtriangledown _{\xi
_{1}}\tciLaplace \xi _{2}-\tciLaplace ^{-1}\bigtriangledown _{\xi
_{2}}\tciLaplace \xi _{1}-\left[ \xi _{1},\xi _{2}\right] \right) \\
&=&\tciLaplace \left( \bigtriangledown _{\xi _{1}}^{\tciLaplace }\xi
_{2}-\bigtriangledown _{\xi _{2}}^{\tciLaplace }\xi _{1}-\left[ \xi _{1},\xi
_{2}\right] \right) \\
&=&\tciLaplace T^{\bigtriangledown ^{\tciLaplace }}\left( \xi _{1},\xi
_{2}\right) .
\end{eqnarray*}

$\left( ii\right) $ The result can be proved the same as in $\left( i\right) 
$.

$\left( iii\right) $ Due to $\left( ii\right) $, it can be easily checked
that%
\begin{eqnarray*}
&&\left( d^{\bigtriangledown }\tciLaplace \right) \left( \xi _{1},\xi
_{2}\right) \\
&=&\left( \bigtriangledown _{\xi _{1}}\tciLaplace \right) \xi _{2}-\left(
\bigtriangledown _{\xi _{2}}\tciLaplace \right) \xi _{1}+\tciLaplace
T^{\bigtriangledown }\left( \xi _{1},\xi _{2}\right) \\
&=&\left( \bigtriangledown _{\xi _{1}}\tciLaplace \right) \xi _{2}-\left(
\bigtriangledown _{\xi _{2}}\tciLaplace \right) \xi _{1}+\left(
d^{\bigtriangledown ^{\tciLaplace }}\tciLaplace \right) \left( \xi _{1},\xi
_{2}\right) .
\end{eqnarray*}%
It is straightforward to obtain%
\begin{equation*}
\left( d^{\bigtriangledown }\tciLaplace \right) \left( \xi _{1},\xi
_{2}\right) -\left( d^{\bigtriangledown ^{\tciLaplace }}\tciLaplace \right)
\left( \xi _{1},\xi _{2}\right) =\left( \bigtriangledown _{\xi
_{1}}\tciLaplace \right) \xi _{2}-\left( \bigtriangledown _{\xi
_{2}}\tciLaplace \right) \xi _{1}
\end{equation*}%
for any vector fields $\xi _{1},\xi _{2}$, that is, $d^{\bigtriangledown
}\tciLaplace =d^{\bigtriangledown ^{\tciLaplace }}\tciLaplace
\Leftrightarrow \left( \bigtriangledown ,\tciLaplace \right) $ is Codazzi
coupled.
\end{proof}

Fei and Zhang \cite{Fei} showed that if a quadric operator $\tciLaplace $
and a linear connection $\nabla $ without torsion are Codazzi-coupled, then
the quadric operator $\tciLaplace $ is integrable. For a linear connection $%
\nabla $ with torsion tensor $T^{\nabla }$, it is possible to give the
following lemma.

\begin{lemma}
\label{lem1}Let $\bigtriangledown $ be a linear connection with torsion
tensor $T^{\nabla }$ and $\tciLaplace $ be a $(1,1)$-tensor field on $\acute{%
N}$. If $d^{\bigtriangledown }\tciLaplace =0$, then the Nijenhuis tensor $%
N_{\tciLaplace }$ associated with $\tciLaplace $ is 
\begin{equation*}
N_{\tciLaplace }\left( \xi _{1},\xi _{2}\right) =-\tciLaplace \left(
T^{\bigtriangledown }\left( \xi _{1},\tciLaplace \xi _{2}\right)
+T^{\bigtriangledown }\left( \tciLaplace \xi _{1},\xi _{2}\right) \right) .
\end{equation*}
\end{lemma}

\begin{proof}
Since $d^{\bigtriangledown }\tciLaplace =0$, there exists the equality 
\begin{equation}
\left[ \xi _{1},\xi _{2}\right] =\tciLaplace ^{-1}\left( \bigtriangledown
_{\xi _{1}}\tciLaplace \xi _{2}-\bigtriangledown _{\xi _{2}}\tciLaplace \xi
_{1}\right) \text{.}  \label{1}
\end{equation}%
From the definition of Nijenhuis tensor $N_{\tciLaplace }$ associated with $%
\tciLaplace $, using the equality (\ref{1}), we can compute%
\begin{eqnarray*}
N_{\tciLaplace }\left( \xi _{1},\xi _{2}\right) &=&-\tciLaplace ^{2}\left[
\xi _{1},\xi _{2}\right] +\tciLaplace \left[ \xi _{1},\tciLaplace \xi _{2}%
\right] +\tciLaplace \left[ \tciLaplace \xi _{1},\xi _{2}\right] -\left[
\tciLaplace \xi _{1},\tciLaplace \xi _{2}\right] \\
&=&-\tciLaplace ^{2}\left( \tciLaplace ^{-1}\left( \bigtriangledown _{\xi
_{1}}\tciLaplace \xi _{2}-\bigtriangledown _{\xi _{2}}\tciLaplace \xi
_{1}\right) \right) +\left( \bigtriangledown _{\xi _{1}}\tciLaplace ^{2}\xi
_{2}-\bigtriangledown _{\tciLaplace \xi _{2}}\tciLaplace \xi _{1}\right) \\
&&+\left( \bigtriangledown _{\tciLaplace \xi _{1}}\tciLaplace \xi
_{2}-\bigtriangledown _{\xi _{2}}\tciLaplace ^{2}\xi _{1}\right)
-\tciLaplace ^{-1}\left( \bigtriangledown _{\tciLaplace \xi _{1}}\tciLaplace
^{2}\xi _{2}-\bigtriangledown _{\tciLaplace \xi _{2}}\tciLaplace ^{2}\xi
_{1}\right) \\
&=&-\tciLaplace \left( \bigtriangledown _{\xi _{1}}\tciLaplace \right) \xi
_{2}+\tciLaplace \left( \bigtriangledown _{\xi _{2}}\tciLaplace \right) \xi
_{1}-\left( \bigtriangledown _{\tciLaplace \xi _{2}}\tciLaplace \right) \xi
_{1}+\left( \bigtriangledown _{\tciLaplace \xi _{1}}\tciLaplace \right) \xi
_{2} \\
&=&\tciLaplace ^{2}T^{\bigtriangledown }\left( \xi _{1},\xi _{2}\right)
-\tciLaplace \left( T^{\bigtriangledown }\left( \xi _{1},\tciLaplace \xi
_{2}\right) +T^{\bigtriangledown }\left( \tciLaplace \xi _{1},\xi
_{2}\right) \right) \\
&&+\tciLaplace \left( \left( \bigtriangledown _{\xi _{1}}\tciLaplace \right)
\xi _{2}-\left( \bigtriangledown _{\xi _{2}}\tciLaplace \right) \xi
_{1}\right) \\
&=&-\tciLaplace \left( T^{\bigtriangledown }\left( \xi _{1},\tciLaplace \xi
_{2}\right) +T^{\bigtriangledown }\left( \tciLaplace \xi _{1},\xi
_{2}\right) \right) .
\end{eqnarray*}
\end{proof}

Consider the condition $T^{\bigtriangledown }\left( \xi _{1},\tciLaplace \xi
_{2}\right) =-T^{\bigtriangledown }\left( \tciLaplace \xi _{1},\xi
_{2}\right) $ which may be called torsion-compatibility. Then we have the
following result.

\begin{proposition}
\label{pro2}An almost complex structure $\tciLaplace $ is integrable if $%
d^{\bigtriangledown }\tciLaplace =0$ and $T^{\bigtriangledown }\left(
\tciLaplace \xi _{1},\xi _{2}\right) =-T^{\bigtriangledown }\left( \xi
_{1},\tciLaplace \xi _{2}\right) $ (torsion-compatibility condition).
\end{proposition}

There is another way to understand the relationship between $%
d^{\bigtriangledown }\tciLaplace $ and integrability of the structure $%
\tciLaplace $. Using the definition of $d^{\bigtriangledown }\tciLaplace $,
it is possible to write the following equality%
\begin{equation*}
\left( d^{\bigtriangledown }\tciLaplace \right) \left( \tciLaplace \xi
_{1},\xi _{2}\right) +\left( d^{\bigtriangledown }\tciLaplace \right) \left(
\xi _{1},\tciLaplace \xi _{2}\right) =T^{\bigtriangledown }\left(
\tciLaplace \xi _{1},\tciLaplace \xi _{2}\right) -T^{\bigtriangledown
}\left( \xi _{1},\xi _{2}\right) -N_{\tciLaplace }\left( \xi _{1},\xi
_{2}\right) .
\end{equation*}%
If the almost complex structure $\tciLaplace $ is integrable, then there
exists the equality $\left( d^{\bigtriangledown }\tciLaplace \right) \left(
\tciLaplace \xi _{1},\xi _{2}\right) =-\left( d^{\bigtriangledown
}\tciLaplace \right) \left( \xi _{1},\tciLaplace \xi _{2}\right) $ provided
that $T^{\bigtriangledown }\left( \tciLaplace \xi _{1},\tciLaplace \xi
_{2}\right) =T^{\bigtriangledown }\left( \xi _{1},\xi _{2}\right) $. Also,
via Proposition \ref{GAD1}, we can give the following results.

\begin{corollary}
Assume that the torsion tensor $T^{\bigtriangledown }$ of a linear
connection $\bigtriangledown $ satisfies the torsion-compatibility
condition. If $T^{\bigtriangledown ^{\tciLaplace }}=0$, then the almost
complex structure $\tciLaplace $ is integrable.
\end{corollary}

It is possible to give an alternative conclusion related to $%
d^{\bigtriangledown }\tciLaplace $ and $T^{\bigtriangledown }$. From the
definition of $d^{\bigtriangledown }\tciLaplace $, we get 
\begin{equation*}
\left( d^{\bigtriangledown }\tciLaplace \right) \left( \tciLaplace \xi
_{1},\xi _{2}\right) =\left( \bigtriangledown _{\tciLaplace \xi
_{1}}\tciLaplace \right) \xi _{2}-\left( \bigtriangledown _{\xi
_{2}}\tciLaplace \right) \tciLaplace \xi _{1}+\tciLaplace
T^{\bigtriangledown }\left( \tciLaplace \xi _{1},\xi _{2}\right)
\end{equation*}%
and%
\begin{equation*}
\left( d^{\bigtriangledown }\tciLaplace \right) \left( \xi _{1},\tciLaplace
\xi _{2}\right) =\left( \bigtriangledown _{\xi _{1}}\tciLaplace \right)
\tciLaplace \xi _{2}-\left( \bigtriangledown _{\tciLaplace \xi
_{2}}\tciLaplace \right) \xi _{1}+\tciLaplace T^{\bigtriangledown }\left(
\xi _{1},\tciLaplace \xi _{2}\right) .
\end{equation*}%
If $\Psi _{\tciLaplace \xi _{1}}\xi _{2}=\nabla _{\tciLaplace \xi _{1}}\xi
_{2}-\tciLaplace \left( \nabla _{\xi _{1}}\xi _{2}\right) =0$ for any vector
fields $\xi _{1}$ and $\xi _{2}$, where $\Psi $ is the Vishnevskii operator 
\cite{SalimovK}, we have%
\begin{equation*}
\left( d^{\bigtriangledown }\tciLaplace \right) \left( \tciLaplace \xi
_{1},\xi _{2}\right) +\left( d^{\bigtriangledown }\tciLaplace \right) \left(
\xi _{1},\tciLaplace \xi _{2}\right) =\tciLaplace \left( T^{\bigtriangledown
}\left( \tciLaplace \xi _{1},\xi _{2}\right) +T^{\bigtriangledown }\left(
\xi _{1},\tciLaplace \xi _{2}\right) \right) ,
\end{equation*}%
from which we immediately obtain 
\begin{equation*}
\left( d^{\bigtriangledown }\tciLaplace \right) \left( \tciLaplace \xi
_{1},\xi _{2}\right) =-\left( d^{\bigtriangledown }\tciLaplace \right)
\left( \xi _{1},\tciLaplace \xi _{2}\right) \text{ }
\end{equation*}%
if and only if the torsion tensor $T^{\bigtriangledown }$ satisfies the
torsion-compatibility condition.

\section{Quasi-statistical structures with a Hermitian metric $g$}

In this section, firstly, we will study quasi-statistical structures
admitting a linear connection $\nabla $ with torsion tensor $T^{\nabla }$, a
pseudo-Riemannian metric $g$ and an almost complex structure $\tciLaplace $,
and obtain interesting results concerning with them. Also, for the K\"{a}%
hler manifolds, we will give a new alternative classification.

\begin{definition}
Let $\acute{N}$ be a differentiable manifold with an almost complex
structure $\tciLaplace $. A Hermitian metric on $\acute{N}$ is a
pseudo-Riemannian metric $g$ such that 
\begin{equation*}
g\left( \tciLaplace \xi _{1},\tciLaplace \xi _{2}\right) =g\left( \xi
_{1},\xi _{2}\right)
\end{equation*}%
or equivalently 
\begin{equation}
g\left( \tciLaplace \xi _{1},\xi _{2}\right) =-g\left( \xi _{1},\tciLaplace
\xi _{2}\right)  \label{GAD0}
\end{equation}%
for any vector fields $\xi _{1}$ and $\xi _{2}$ on $\acute{N}$. Then the
triple $\left( \acute{N},g,\tciLaplace \right) $ is an almost Hermitian
manifold. The fundamental $2$-form $\omega $ is given by $\omega \left( \xi
_{1},\xi _{2}\right) =g\left( \tciLaplace \xi _{1},\xi _{2}\right) $ for any
vector fields $\xi _{1},\xi _{2}$ on $\acute{N}$. $\left( g,\tciLaplace
,\omega \right) $ is known as the \textquotedblleft compatible
triple.\textquotedblright\ If the almost complex structure $\tciLaplace $ is
integrable, the triple $\left( \acute{N},g,\tciLaplace \right) $ is a
Hermitian manifold. The triple $\left( \acute{N},\omega ,\tciLaplace \right) 
$ is a K\"{a}hler manifold if the structure $\tciLaplace $ is integrable and 
$\omega $ is closed, that is, $d\omega =0$ or equivalent to these two
conditions is that the structure $\tciLaplace $ is covariantly constant with
respect to the Levi-Civita connection $\bigtriangledown ^{g}$ of $g$ \cite%
{Gray}.
\end{definition}

The well-known formula of the covariant derivative of $g$ with respect to $%
\nabla $ is as follow%
\begin{equation*}
\left( \bigtriangledown _{\xi _{3}}g\right) \left( \xi _{1},\xi _{2}\right)
=\xi _{3}g\left( \xi _{1},\xi _{2}\right) -g\left( \bigtriangledown _{\xi
_{3}}\xi _{1},\xi _{2}\right) -g\left( \xi _{1},\bigtriangledown _{\xi
_{3}}\xi _{2}\right) .
\end{equation*}%
Clearly $\left( \bigtriangledown _{\xi _{3}}g\right) \left( \xi _{1},\xi
_{2}\right) =\left( \bigtriangledown _{\xi _{3}}g\right) \left( \xi _{2},\xi
_{1}\right) ,$ due to symmetry of $g$. It is clear that $g$ is parallel
under $\bigtriangledown $ if and only if $\bigtriangledown g=0$.\newline

Given a pair $\left( \bigtriangledown ,g\right) $, we can also construct $%
\bigtriangledown ^{\ast }$, called a $g$-conjugate connection by%
\begin{equation*}
\xi _{3}g\left( \xi _{1},\xi _{2}\right) =g\left( \bigtriangledown _{\xi
_{3}}\xi _{1},\xi _{2}\right) +g\left( \xi _{1},\bigtriangledown _{\xi
_{3}}^{\ast }\xi _{2}\right) .
\end{equation*}%
It is easy to see that $\bigtriangledown ^{\ast }$ is a linear connection
and a $g$-conjugate of a connection $\nabla $ is involutive, that is, $%
\left( \bigtriangledown ^{\ast }\right) ^{\ast }=\bigtriangledown $. These
two constructions from an arbitrary pair $\left( \bigtriangledown ,g\right) $
are related via $\left( \bigtriangledown _{\xi _{3}}g\right) \left( \xi
_{1},\xi _{2}\right) =g\left( \left( \bigtriangledown ^{\ast
}-\bigtriangledown \right) _{\xi _{3}}\xi _{1},\xi _{2}\right) $, which
satisfy%
\begin{equation*}
\left( \bigtriangledown _{\xi _{3}}^{\ast }g\right) \left( \xi _{1},\xi
_{2}\right) =-\left( \bigtriangledown _{\xi _{3}}g\right) \left( \xi
_{1},\xi _{2}\right) .
\end{equation*}%
Therefore, we say that $\left( \bigtriangledown _{\xi _{3}}^{\ast }g\right)
\left( \xi _{1},\xi _{2}\right) =\left( \bigtriangledown _{\xi _{3}}g\right)
\left( \xi _{1},\xi _{2}\right) =0$ if and only if $\bigtriangledown ^{\ast
}=\bigtriangledown $, that is, $\bigtriangledown $ is $g$-self conjugate. A
linear connection that is both $g$-self conjugate and torsion free is the
Levi Civita connection $\bigtriangledown ^{g}$ of $g$.

\begin{definition}
Let $\bigtriangledown $ be a torsion free linear connection on the
pseudo-Riemannian manifold $\left( \acute{N},g\right) $ with a
pseudo-Riemannian metric $g$. We can say that $\left( \acute{N}%
,g,\bigtriangledown \right) $ is a statistical manifold if the following
equation is satisfied \cite{Lauritzen} 
\begin{equation*}
\left( \bigtriangledown _{\xi _{1}}g\right) \left( \xi _{2},\xi _{3}\right)
=\left( \bigtriangledown _{\xi _{2}}g\right) \left( \xi _{1},\xi _{3}\right)
.
\end{equation*}%
We will consider an extension of the notion of a statistical structure. We
can say that $\left( \acute{N},g,\bigtriangledown \right) $ is a statistical
manifold admitting torsion (SMAT) if $d^{\bigtriangledown }g=0$, where 
\begin{equation*}
\left( d^{\bigtriangledown }g\right) \left( \xi _{1},\xi _{2},\xi
_{3}\right) =\left( \bigtriangledown _{\xi _{1}}g\right) \left( \xi _{2},\xi
_{3}\right) -\left( \bigtriangledown _{\xi _{2}}g\right) \left( \xi _{1},\xi
_{3}\right) +g\left( T^{\bigtriangledown }\left( \xi _{1},\xi _{2}\right)
,\xi _{3}\right)
\end{equation*}%
for any vector fields $\xi _{1},\xi _{2}$ and $\xi _{3}$. Also, it is called
a statistical manifold admitting torsion (SMAT) as a quasi-statistical
manifold \cite{Furuhata}.
\end{definition}

The fundamental $2$-form $\omega $ on $\acute{N}$ is also an almost
symplectic structure. Let us introduce the $\omega $-conjugate
transformation $\bigtriangledown ^{\dag }$ of $\bigtriangledown $ by 
\begin{equation*}
\xi _{3}\omega \left( \xi _{1},\xi _{2}\right) =\omega \left(
\bigtriangledown _{\xi _{3}}\xi _{1},\xi _{2}\right) +\omega \left( \xi
_{1},\bigtriangledown _{\xi _{3}}^{\dag }\xi _{2}\right) ,
\end{equation*}%
where conjugation is invariantly defined with respect to either the first or
the second entry of $\omega $ despite of the skew-symmetric nature of $%
\omega $ \cite{Fei} . The covariant derivative of $\omega $ with respect to $%
\nabla $ is the following $\left( 0,3\right) $-tensor field%
\begin{equation*}
\left( \bigtriangledown _{\xi _{3}}\omega \right) \left( \xi _{1},\xi
_{2}\right) =\xi _{3}\omega \left( \xi _{1},\xi _{2}\right) -\omega \left(
\bigtriangledown _{\xi _{3}}\xi _{1},\xi _{2}\right) -\omega \left( \xi
_{1},\bigtriangledown _{\xi _{3}}\xi _{2}\right) ,
\end{equation*}%
which is skew-symmetric in $\xi _{1},\xi _{2}:$ $\left( \bigtriangledown
_{\xi _{3}}\omega \right) \left( \xi _{1},\xi _{2}\right) =-\left(
\bigtriangledown _{\xi _{3}}\omega \right) \left( \xi _{2},\xi _{1}\right) $%
. Imposing the Codazzi coupling condition of $\left( \bigtriangledown
,\omega \right) $, that is, $\left( \bigtriangledown _{\xi _{3}}\omega
\right) \left( \xi _{1},\xi _{2}\right) =\left( \bigtriangledown _{\xi
_{1}}\omega \right) \left( \xi _{3},\xi _{2}\right) $ leads to $\left(
\bigtriangledown _{\xi _{3}}\omega \right) \left( \xi _{1},\xi _{2}\right)
=0 $.

\begin{lemma}
\label{lem2} \cite{Grigorian} Let $\left( \acute{N},\omega \right) $ be an
almost symplectic manifold with the fundamental $2$-form $\omega $. Then,
for any vector fields $\xi _{1},\xi _{2}$ and $\xi _{3}$,%
\begin{eqnarray*}
d\omega \left( \xi _{1},\xi _{2},\xi _{3}\right) &=&\left( \bigtriangledown
_{\xi _{3}}\omega \right) \left( \xi _{1},\xi _{2}\right) +\left(
\bigtriangledown _{\xi _{1}}\omega \right) \left( \xi _{2},\xi _{3}\right)
+\left( \bigtriangledown _{\xi _{2}}\omega \right) \left( \xi _{3},\xi
_{1}\right) \\
&&+\omega \left( T^{\bigtriangledown }\left( \xi _{1},\xi _{2}\right) ,\xi
_{3}\right) +\omega \left( T^{\bigtriangledown }\left( \xi _{2},\xi
_{3}\right) ,\xi _{1}\right) +\omega \left( T^{\bigtriangledown }\left( \xi
_{3},\xi _{1}\right) ,\xi _{2}\right) .
\end{eqnarray*}
\end{lemma}

\begin{proposition}
\label{pro3}Let $\left( \acute{N},g\right) $ be a pseudo-Riemannian manifold
and let $\bigtriangledown $ be a linear connection with torsion tensor $%
T^{\nabla }$. Let $\omega $ be the fundamental $2$-form on $\acute{N}$.
Then, there exist the following expressions

$\left( i\right) $ Assume that $\left( \bigtriangledown ,\tciLaplace \right) 
$ is Codazzi-coupled. $d^{\bigtriangledown ^{\ast }}\omega =0\Leftrightarrow
\left( \bigtriangledown ^{\ast },g\right) $ is a quasi-statistical structure.

$\left( ii\right) $ Assume that $\left( \bigtriangledown ,\tciLaplace
\right) $ is Codazzi-coupled. $d^{\bigtriangledown ^{\dag }}\omega
=0\Leftrightarrow \left( \bigtriangledown ^{\dag },g\right) $ is a
quasi-statistical structure.

$\left( iii\right) $ Assume that $\left( \bigtriangledown ^{\ast
},\tciLaplace \right) $ is Codazzi-coupled. $d^{\bigtriangledown }\omega
=0\Leftrightarrow \left( \bigtriangledown ,g\right) $ is a quasi-statistical
structure.

$\left( iv\right) $ Assume that $\left( \bigtriangledown ^{\dag
},\tciLaplace \right) $ is Codazzi-coupled. $d^{\bigtriangledown }\omega
=0\Leftrightarrow \left( \bigtriangledown ,g\right) $ is a quasi-statistical
structure.

$\left( v\right) $ Assume that $\left( \bigtriangledown ^{\dag },\tciLaplace
\right) $ is Codazzi-coupled. $d^{\bigtriangledown ^{\tciLaplace }}\omega
=0\Leftrightarrow d^{\bigtriangledown }\omega =0$.

$\left( vi\right) $ Assume that $\left( \bigtriangledown ^{\ast
},\tciLaplace \right) $ is Codazzi-coupled. $\left( \bigtriangledown
^{\tciLaplace },g\right) $ is a quasi-statistical structure $\Leftrightarrow
\left( \bigtriangledown ,g\right) $ is a quasi-statistical structure.
\end{proposition}

\begin{proof}
$\left( i\right) $ We can write%
\begin{eqnarray*}
&&\left( d^{\bigtriangledown ^{\ast }}\omega \right) \left( \xi _{1},\xi
_{2},\xi _{3}\right) \\
&=&\left( \bigtriangledown _{\xi _{1}}^{\ast }\omega \right) \left( \xi
_{2},\xi _{3}\right) -\left( \bigtriangledown _{\xi _{2}}^{\ast }\omega
\right) \left( \xi _{1},\xi _{3}\right) +\omega \left( T^{\bigtriangledown
^{\ast }}\left( \xi _{1},\xi _{2}\right) ,\xi _{3}\right) \\
&=&\xi _{1}\omega \left( \xi _{2},\xi _{3}\right) -\omega \left(
\bigtriangledown _{\xi _{1}}^{\ast }\xi _{2},\xi _{3}\right) -\omega \left(
\xi _{2},\bigtriangledown _{\xi _{1}}^{\ast }\xi _{3}\right) -\xi _{2}\omega
\left( \xi _{1},\xi _{3}\right) \\
&&+\omega \left( \bigtriangledown _{\xi _{2}}^{\ast }\xi _{1},\xi
_{3}\right) +\omega \left( \xi _{1},\bigtriangledown _{\xi _{2}}^{\ast }\xi
_{3}\right) +\omega \left( T^{\bigtriangledown ^{\ast }}\left( \xi _{1},\xi
_{2}\right) ,\xi _{3}\right) \\
&=&\xi _{1}g\left( \tciLaplace \xi _{2},\xi _{3}\right) -g\left( \tciLaplace
\bigtriangledown _{\xi _{1}}^{\ast }\xi _{2},\xi _{3}\right) -g\left(
\tciLaplace \xi _{2},\bigtriangledown _{\xi _{1}}^{\ast }\xi _{3}\right)
-\xi _{2}g\left( \tciLaplace \xi _{1},\xi _{3}\right) \\
&&+g\left( \tciLaplace \bigtriangledown _{\xi _{2}}^{\ast }\xi _{1},\xi
_{3}\right) +g\left( \tciLaplace \xi _{1},\bigtriangledown _{\xi _{2}}^{\ast
}\xi _{3}\right) +g\left( \tciLaplace T^{\bigtriangledown ^{\ast }}\left(
\xi _{1},\xi _{2}\right) ,\xi _{3}\right) \\
&=&-\xi _{1}g\left( \xi _{2},\tciLaplace \xi _{3}\right) +g\left(
\bigtriangledown _{\xi _{1}}^{\ast }\xi _{2},\tciLaplace \xi _{3}\right)
+g\left( \xi _{2},\tciLaplace \bigtriangledown _{\xi _{1}}^{\ast }\xi
_{3}\right) +\xi _{2}g\left( \xi _{1},\tciLaplace \xi _{3}\right) \\
&&-g\left( \bigtriangledown _{\xi _{2}}^{\ast }\xi _{1},\tciLaplace \xi
_{3}\right) -g\left( \xi _{1},\tciLaplace \bigtriangledown _{\xi _{2}}^{\ast
}\xi _{3}\right) -g\left( T^{\bigtriangledown ^{\ast }}\left( \xi _{1},\xi
_{2}\right) ,\tciLaplace \xi _{3}\right) \\
&=&-\left( \bigtriangledown _{\xi _{1}}^{\ast }g\right) \left( \xi
_{2},\tciLaplace \xi _{3}\right) +\left( \bigtriangledown _{\xi _{2}}^{\ast
}g\right) \left( \xi _{1},\tciLaplace \xi _{3}\right) -g\left(
T^{\bigtriangledown ^{\ast }}\left( \xi _{1},\xi _{2}\right) ,\tciLaplace
\xi _{3}\right) \\
&&-g\left( \xi _{2},\left( \bigtriangledown _{\xi _{1}}^{\ast }\tciLaplace
\right) \xi _{3}\right) +g\left( \xi _{1},\left( \bigtriangledown _{\xi
_{2}}^{\ast }\tciLaplace \right) \xi _{3}\right) \\
&=&-\left( d^{\bigtriangledown ^{\ast }}g\right) \left( \xi _{1},\xi
_{2},\tciLaplace \xi _{3}\right) -g\left( \xi _{2},\left( \bigtriangledown
_{\xi _{1}}^{\ast }\tciLaplace \right) \xi _{3}\right) +g\left( \xi
_{1},\left( \bigtriangledown _{\xi _{2}}^{\ast }\tciLaplace \right) \xi
_{3}\right) .
\end{eqnarray*}%
Since $g\left( \xi _{1},\left( \bigtriangledown _{\xi _{2}}^{\ast
}\tciLaplace \right) \xi _{3}\right) =-g\left( \xi _{3},\left(
\bigtriangledown _{\xi _{2}}\tciLaplace \right) \xi _{1}\right) $ and the
pair $\left( \bigtriangledown ,\tciLaplace \right) $ is Codazzi-coupled, we
get 
\begin{equation*}
\left( d^{\bigtriangledown ^{\ast }}\omega \right) \left( \xi _{1},\xi
_{2},\xi _{3}\right) =-\left( d^{\bigtriangledown ^{\ast }}g\right) \left(
\xi _{1},\xi _{2},\tciLaplace \xi _{3}\right) ,
\end{equation*}%
that is, $d^{\bigtriangledown ^{\ast }}\omega =0\Leftrightarrow \left(
\bigtriangledown ^{\ast },g\right) $ is a quasi-statistical structure.

$\left( ii\right) $ We obtain the following%
\begin{eqnarray*}
&&\left( d^{\bigtriangledown ^{\dag }}\omega \right) \left( \xi _{1},\xi
_{2},\xi _{3}\right) \\
&=&\left( \bigtriangledown _{\xi _{1}}^{\dag }\omega \right) \left( \xi
_{2},\xi _{3}\right) -\left( \bigtriangledown _{\xi _{2}}^{\dag }\omega
\right) \left( \xi _{1},\xi _{3}\right) +\omega \left( T^{\bigtriangledown
^{\dag }}\left( \xi _{1},\xi _{2}\right) ,\xi _{3}\right) \\
&=&\xi _{1}\omega \left( \xi _{2},\xi _{3}\right) -\omega \left( \xi
_{2},\bigtriangledown _{\xi _{1}}^{\dag }\xi _{3}\right) -\xi _{2}\omega
\left( \xi _{1},\xi _{3}\right) +\omega \left( \xi _{1},\bigtriangledown
_{\xi _{2}}^{\dag }\xi _{3}\right) \\
&&-\omega \left( \left[ \xi _{1},\xi _{2}\right] ,\xi _{3}\right) \\
&=&\xi _{1}g\left( \tciLaplace \xi _{2},\xi _{3}\right) -g\left( \tciLaplace
\xi _{2},\bigtriangledown _{\xi _{1}}^{\dag }\xi _{3}\right) -\xi
_{2}g\left( \tciLaplace \xi _{1},\xi _{3}\right) +g\left( \tciLaplace \xi
_{1},\bigtriangledown _{\xi _{2}}^{\dag }\xi _{3}\right) \\
&&-g\left( \tciLaplace \left[ \xi _{1},\xi _{2}\right] ,\xi _{3}\right) \\
&=&-\xi _{1}g\left( \xi _{2},\tciLaplace \xi _{3}\right) +g\left( \xi
_{2},\tciLaplace \bigtriangledown _{\xi _{1}}^{\dag }\xi _{3}\right) +\xi
_{2}g\left( \xi _{1},\tciLaplace \xi _{3}\right) -g\left( \xi
_{1},\tciLaplace \bigtriangledown _{\xi _{2}}^{\dag }\xi _{3}\right) \\
&&+g\left( \left[ \xi _{1},\xi _{2}\right] ,\tciLaplace \xi _{3}\right) \\
&=&-\left( \bigtriangledown _{\xi _{1}}^{\dag }g\right) \left( \xi
_{2},\tciLaplace \xi _{3}\right) +\left( \bigtriangledown _{\xi _{2}}^{\dag
}g\right) \left( \xi _{1},\tciLaplace \xi _{3}\right) -g\left(
T^{\bigtriangledown ^{\dag }}\left( \xi _{1},\xi _{2}\right) ,\tciLaplace
\xi _{3}\right) \\
&&-g\left( \xi _{2},\left( \bigtriangledown _{\xi _{1}}^{\dag }\tciLaplace
\right) \xi _{3}\right) +g\left( \xi _{1},\left( \bigtriangledown _{\xi
_{2}}^{\dag }\tciLaplace \right) \xi _{3}\right) \\
&=&-\left( d^{\bigtriangledown ^{\dag }}g\right) \left( \xi _{1},\xi
_{2},\tciLaplace \xi _{3}\right) -g\left( \xi _{2},\left( \bigtriangledown
_{\xi _{1}}^{\dag }\tciLaplace \right) \xi _{3}\right) +g\left( \xi
_{1},\left( \bigtriangledown _{\xi _{2}}^{\dag }\tciLaplace \right) \xi
_{3}\right)
\end{eqnarray*}%
such that, from the hypothesis, we have%
\begin{equation*}
\left( d^{\bigtriangledown \dag }\omega \right) \left( \xi _{1},\xi _{2},\xi
_{3}\right) =-\left( d^{\bigtriangledown ^{\dag }}g\right) \left( \xi
_{1},\xi _{2},\tciLaplace \xi _{3}\right) ,
\end{equation*}%
that is, $d^{\bigtriangledown ^{\dag }}\omega =0\Leftrightarrow \left(
\bigtriangledown ^{\dag },g\right) $ is a quasi-statistical structure.

$\left( iii\right) $ The result can be proved the same as in $\left(
i\right) $.

$\left( iv\right) $ The result can be proved the same as in $\left(
ii\right) $.

$\left( v\right) $ For any vector fields $\xi _{1},\xi _{2}$ and $\xi _{3}$,
we have the following%
\begin{eqnarray*}
&&\left( d^{\bigtriangledown ^{\tciLaplace }}\omega \right) \left( \xi
_{1},\xi _{2},\xi _{3}\right) \\
&=&\left( \bigtriangledown _{\xi _{1}}^{\tciLaplace }\omega \right) \left(
\xi _{2},\xi _{3}\right) -\left( \bigtriangledown _{\xi _{2}}^{\tciLaplace
}\omega \right) \left( \xi _{1},\xi _{3}\right) +\omega \left(
T^{\bigtriangledown ^{\tciLaplace }}\left( \xi _{1},\xi _{2}\right) ,\xi
_{3}\right) \\
&=&\xi _{1}\omega \left( \xi _{2},\xi _{3}\right) -\omega \left( \xi
_{2},\bigtriangledown _{\xi _{1}}^{\tciLaplace }\xi _{3}\right) -\xi
_{2}\omega \left( \xi _{1},\xi _{3}\right) \\
&&+\omega \left( \xi _{1},\bigtriangledown _{\xi _{2}}^{\tciLaplace }\xi
_{3}\right) -\omega \left( \left[ \xi _{1},\xi _{2}\right] ,\xi _{3}\right)
\\
&=&\xi _{1}\omega \left( \xi _{2},\xi _{3}\right) -\omega \left( \xi
_{2},\tciLaplace ^{-1}\bigtriangledown _{\xi _{1}}\tciLaplace \xi
_{3}\right) -\xi _{2}\omega \left( \xi _{1},\xi _{3}\right) \\
&&+\omega \left( \xi _{1},\tciLaplace ^{-1}\bigtriangledown _{\xi
_{2}}\tciLaplace \xi _{3}\right) -\omega \left( \left[ \xi _{1},\xi _{2}%
\right] ,\xi _{3}\right) \\
&=&\left( \bigtriangledown _{\xi _{1}}\omega \right) \left( \xi _{2},\xi
_{3}\right) -\left( \bigtriangledown _{\xi _{2}}\omega \right) \left( \xi
_{1},\xi _{3}\right) +\omega \left( T^{\bigtriangledown }\left( \xi _{1},\xi
_{2}\right) ,\xi _{3}\right) \\
&&-\omega \left( \xi _{2},\tciLaplace ^{-1}\left( \bigtriangledown _{\xi
_{1}}\tciLaplace \right) \xi _{3}\right) +\omega \left( \xi _{1},\tciLaplace
^{-1}\left( \bigtriangledown _{\xi _{2}}\tciLaplace \right) \xi _{3}\right)
\\
&=&\left( d^{\bigtriangledown }\omega \right) \left( \xi _{1},\xi _{2},\xi
_{3}\right) -\omega \left( \xi _{2},\tciLaplace ^{-1}\left( \bigtriangledown
_{\xi _{1}}\tciLaplace \right) \xi _{3}\right) +\omega \left( \xi
_{1},\tciLaplace ^{-1}\left( \bigtriangledown _{\xi _{2}}\tciLaplace \right)
\xi _{3}\right) .
\end{eqnarray*}%
From the equation $\omega \left( \xi _{1},\left( \bigtriangledown _{\xi
_{2}}\tciLaplace \right) \xi _{3}\right) =\omega \left( \xi _{3},\left(
\bigtriangledown _{\xi _{2}}^{\dag }\tciLaplace \right) \xi _{1}\right) $,
it is easy to see that%
\begin{equation*}
\left( d^{\bigtriangledown ^{\tciLaplace }}\omega \right) \left( \xi
_{1},\xi _{2},\xi _{3}\right) =\left( d^{\bigtriangledown }\omega \right)
\left( \xi _{1},\xi _{2},\xi _{3}\right) ,
\end{equation*}%
that is,%
\begin{equation*}
d^{\bigtriangledown ^{\tciLaplace }}\omega =0\Leftrightarrow
d^{\bigtriangledown }\omega =0.
\end{equation*}

$\left( vi\right) $ From the definition of the $\tciLaplace $-conjugate
transformation, we write the following%
\begin{eqnarray*}
&&\left( d^{\bigtriangledown ^{\tciLaplace }}g\right) \left( \xi _{1},\xi
_{2},\xi _{3}\right) \\
&=&\left( \bigtriangledown _{\xi _{1}}^{\tciLaplace }g\right) \left( \xi
_{2},\xi _{3}\right) -\left( \bigtriangledown _{\xi _{2}}^{\tciLaplace
}g\right) \left( \xi _{1},\xi _{3}\right) +g\left( T^{\bigtriangledown
^{\tciLaplace }}\left( \xi _{1},\xi _{2}\right) ,\xi _{3}\right) \\
&=&\xi _{1}g\left( \xi _{2},\xi _{3}\right) -g\left( \xi
_{2},\bigtriangledown _{\xi _{1}}^{\tciLaplace }\xi _{3}\right) -\xi
_{2}g\left( \xi _{1},\xi _{3}\right) \\
&&+g\left( \xi _{1},\bigtriangledown _{\xi _{2}}^{\tciLaplace }\xi
_{3}\right) -g\left( \left[ \xi _{1},\xi _{2}\right] ,\xi _{3}\right) \\
&=&\left( \bigtriangledown _{\xi _{1}}g\right) \left( \xi _{2},\xi
_{3}\right) -\left( \bigtriangledown _{\xi _{2}}g\right) \left( \xi _{1},\xi
_{3}\right) +g\left( T^{\bigtriangledown }\left( \xi _{1},\xi _{2}\right)
,\xi _{3}\right) \\
&&-g\left( \xi _{2},\tciLaplace ^{-1}\left( \bigtriangledown _{\xi
_{1}}\tciLaplace \right) \xi _{3}\right) +g\left( \xi _{1},\tciLaplace
^{-1}\left( \bigtriangledown _{\xi _{2}}\tciLaplace \right) \xi _{3}\right)
\\
&=&\left( d^{\bigtriangledown }g\right) \left( \xi _{1},\xi _{2},\xi
_{3}\right) -g\left( \xi _{2},\tciLaplace ^{-1}\left( \bigtriangledown _{\xi
_{1}}\tciLaplace \right) \xi _{3}\right) +g\left( \xi _{1},\tciLaplace
^{-1}\left( \bigtriangledown _{\xi _{2}}\tciLaplace \right) \xi _{3}\right) .
\end{eqnarray*}%
To complete the proof, one needs to note%
\begin{equation*}
g\left( \xi _{1},\left( \bigtriangledown _{\xi _{2}}^{\ast }\tciLaplace
\right) \xi _{3}\right) =-g\left( \xi _{3},\left( \bigtriangledown _{\xi
_{2}}\tciLaplace \right) \xi _{1}\right) .
\end{equation*}%
Immediately, from the hypothesis, we say that%
\begin{equation*}
\left( d^{\bigtriangledown ^{\tciLaplace }}g\right) \left( \xi _{1},\xi
_{2},\xi _{3}\right) =\left( d^{\bigtriangledown }g\right) \left( \xi
_{1},\xi _{2},\xi _{3}\right) ,
\end{equation*}%
i.e., $\left( \bigtriangledown ^{\tciLaplace },g\right) $ is a
quasi-statistical structure if and only if $\left( \bigtriangledown
,g\right) $ is a quasi-statistical structure.
\end{proof}

\begin{remark}
\label{teo1} Let $\left( \acute{N},g\right) $ be a pseudo-Riemannian
manifold and let $\bigtriangledown $ be a linear connection with torsion
tensor $T^{\nabla }$. Let $\omega $ be the fundamental $2$-form on $\acute{N}
$. Denote by $\bigtriangledown ^{\ast },\bigtriangledown ^{\dag }$ and $%
\bigtriangledown ^{\tciLaplace }$, $g$-conjugate, $\omega $-conjugate and $%
\tciLaplace $-conjugate transformations of a linear connection $%
\bigtriangledown $. These transformations are involutive $\left(
\bigtriangledown ^{\ast }\right) ^{\ast }=\left( \bigtriangledown ^{\dag
}\right) ^{\dag }=\left( \bigtriangledown ^{\tciLaplace }\right)
^{\tciLaplace }=\bigtriangledown $. When the equation (\ref{GAD0}) is
satisfied, these transformations of $\nabla $ are commutative%
\begin{equation*}
\bigtriangledown ^{\ast }=\left( \bigtriangledown ^{\dag }\right)
^{\tciLaplace }=\left( \bigtriangledown ^{\tciLaplace }\right) ^{\dag },
\end{equation*}%
\begin{equation*}
\bigtriangledown ^{\dag }=\left( \bigtriangledown ^{\ast }\right)
^{\tciLaplace }=\left( \bigtriangledown ^{\tciLaplace }\right) ^{\ast },
\end{equation*}%
\begin{equation*}
\bigtriangledown ^{\tciLaplace }=\left( \bigtriangledown ^{\ast }\right)
^{\dag }=\left( \bigtriangledown ^{\dag }\right) ^{\ast }.
\end{equation*}%
Hence, $\left( id,\ast ,\dag ,\tciLaplace \right) $ forms a $4$-element
Klein group of transformation of linear connections on $\acute{N}$ ( for
details, see Theorem 2.13 in \cite{Fei}).
\end{remark}

\begin{proposition}
\label{pro4}Let $\left( \acute{N},g\right) $ be a pseudo-Riemannian manifold
and let $\bigtriangledown $ be a linear connection with torsion tensor $%
T^{\nabla }$ on $\acute{N}$. Let $\omega $ be the fundamental $2$-form on $%
\acute{N}$. Then, the following expressions hold

$\left( i\right) $ $d^{\bigtriangledown ^{\tciLaplace }}\omega
=0\Leftrightarrow \left( \bigtriangledown ,g\right) $ is a quasi-statistical
structure;

$\left( ii\right) $ $d^{\bigtriangledown }\omega =0\Leftrightarrow \left(
\bigtriangledown ^{\tciLaplace },g\right) $ is a quasi-statistical structure;

$\left( iii\right) $ $d^{\bigtriangledown ^{\dag }}\omega =0\Leftrightarrow
\left( \bigtriangledown ^{\ast },g\right) $ is a quasi-statistical structure;

$\left( iv\right) $ $d^{\bigtriangledown ^{\ast }}\omega =0\Leftrightarrow
\left( \bigtriangledown ^{\dag },g\right) $ is a quasi-statistical structure.
\end{proposition}

\begin{proof}
To show $\left( i\right) $, one only needs the following equality%
\begin{eqnarray*}
&&\left( d^{\bigtriangledown ^{\tciLaplace }}\omega \right) \left( \xi
_{1},\xi _{2},\xi _{3}\right) \\
&=&\left( \bigtriangledown _{\xi _{1}}^{\tciLaplace }\omega \right) \left(
\xi _{2},\xi _{3}\right) -\left( \bigtriangledown _{\xi _{2}}^{\tciLaplace
}\omega \right) \left( \xi _{1},\xi _{3}\right) +\omega \left(
T^{\bigtriangledown ^{\tciLaplace }}\left( \xi _{1},\xi _{2}\right) ,\xi
_{3}\right) \\
&=&\xi _{1}\omega \left( \xi _{2},\xi _{3}\right) -\omega \left( \xi
_{2},\bigtriangledown _{\xi _{1}}^{\tciLaplace }\xi _{3}\right) -\xi
_{2}\omega \left( \xi _{1},\xi _{3}\right) \\
&&+\omega \left( \xi _{1},\bigtriangledown _{\xi _{2}}^{\tciLaplace }\xi
_{3}\right) -\omega \left( \left[ \xi _{1},\xi _{2}\right] ,\xi _{3}\right)
\\
&=&\xi _{1}\omega \left( \xi _{2},\xi _{3}\right) -\omega \left( \xi
_{2},\tciLaplace ^{-1}\bigtriangledown _{\xi _{1}}\tciLaplace \xi
_{3}\right) -\xi _{2}\omega \left( \xi _{1},\xi _{3}\right) \\
&&+\omega \left( \xi _{1},\tciLaplace ^{-1}\bigtriangledown _{\xi
_{2}}\tciLaplace \xi _{3}\right) -\omega \left( \left[ \xi _{1},\xi _{2}%
\right] ,\xi _{3}\right) \\
&=&\xi _{1}g\left( \tciLaplace \xi _{2},\xi _{3}\right) -g\left( \tciLaplace
\xi _{2},\tciLaplace ^{-1}\bigtriangledown _{\xi _{1}}\tciLaplace \xi
_{3}\right) -\xi _{2}g\left( \tciLaplace \xi _{1},\xi _{3}\right) \\
&&+g\left( \tciLaplace \xi _{1},\tciLaplace ^{-1}\bigtriangledown _{\xi
_{2}}\tciLaplace \xi _{3}\right) -g\left( \tciLaplace \left[ \xi _{1},\xi
_{2}\right] ,\xi _{3}\right) \\
&=&-\xi _{1}g\left( \xi _{2},\tciLaplace \xi _{3}\right) +g\left( \xi
_{2},\bigtriangledown _{\xi _{1}}\tciLaplace \xi _{3}\right) +\xi
_{2}g\left( \xi _{1},\tciLaplace \xi _{3}\right) \\
&&-g\left( \xi _{1},\bigtriangledown _{\xi _{2}}\tciLaplace \xi _{3}\right)
+g\left( \left[ \xi _{1},\xi _{2}\right] ,\tciLaplace \xi _{3}\right) \\
&=&-\left( \bigtriangledown _{\xi _{1}}g\right) \left( \xi _{2},\tciLaplace
\xi _{3}\right) +\left( \bigtriangledown _{\xi _{2}}g\right) \left( \xi
_{1},\tciLaplace \xi _{3}\right) -g\left( T^{\bigtriangledown }\left( \xi
_{1},\xi _{2}\right) ,\tciLaplace \xi _{3}\right) \\
&=&-\left( d^{\bigtriangledown }g\right) \left( \xi _{1},\xi
_{2},\tciLaplace \xi _{3}\right) ,
\end{eqnarray*}%
which completes the proof, i.e., $d^{\bigtriangledown ^{\tciLaplace }}\omega
=0\Leftrightarrow \left( \bigtriangledown ,g\right) $ is a quasi-statistical
structure. From Remark \ref{teo1}, the other statements can easily be proved.
\end{proof}

\begin{corollary}
\label{cor4}Let $\acute{N}$ be a pseudo-Riemannian manifold equipped with a
pseudo Riemannian metric $g$ and a linear connection $\bigtriangledown $.
Let $\left( g,\omega ,\tciLaplace \right) $ be a compatible triple, and $%
\bigtriangledown ^{\ast },\bigtriangledown ^{\dag }$and $\bigtriangledown
^{\tciLaplace }$ denote, respectively, $g$-conjugate, $\omega $-conjugate
and $\tciLaplace $-conjugate transformations of an arbitrary linear
connection $\bigtriangledown $. Fei and Zhang \cite{Fei} showed that $\left(
\bigtriangledown ,g\right) $ is a statistical structure if and only if $%
\bigtriangledown ^{\ast }$ is torsion-free. There exist the following
expressions

$\left( i\right) $ $\left( \bigtriangledown ^{\ast },g\right) $ is a
quasi-statistical structure if and only if $\bigtriangledown $ is
torsion-free;

$\left( ii\right) $ $\left( \bigtriangledown ^{\tciLaplace },g\right) $ is a
quasi-statistical structure if and only if $\bigtriangledown ^{\dag }$ is
torsion-free;

$\left( iii\right) $ $\left( \bigtriangledown ^{\dag },g\right) $ is a
quasi-statistical structure if and only if $\bigtriangledown ^{\tciLaplace }$
is torsion-free.
\end{corollary}

\begin{proof}
From $4$-element Klein group action $\left( id,\ast ,\dag ,\tciLaplace
\right) $ on the space of linear connections, the results immediately follow.
\end{proof}

As a corollary of Proposition \ref{pro4} and Corollary \ref{cor4}, we get
the following.

\begin{proposition}
\label{pro5} Let $\acute{N}$ be a pseudo-Riemannian manifold equipped with a
pseudo-Riemannian metric $g$ and a linear connection $\bigtriangledown $.
Let $\left( g,\omega ,\tciLaplace \right) $ be a compatible triple, and $%
\bigtriangledown ^{\ast },\bigtriangledown ^{\dag }$and $\bigtriangledown
^{\tciLaplace }$ denote, respectively, $g$-conjugate, $\omega $-conjugate
and $\tciLaplace $-conjugate transformations of the linear connection $%
\bigtriangledown $. Then, there exist the followings

$\left( i\right) $ $d^{\bigtriangledown }\omega =0\Leftrightarrow
T^{\bigtriangledown ^{\dag }}=0\Leftrightarrow d^{\bigtriangledown ^{\ast
}}\tciLaplace =0\Leftrightarrow d^{\bigtriangledown ^{\tciLaplace }}g=0$;

$\left( ii\right) $ $d^{\bigtriangledown ^{\ast }}\omega =0\Leftrightarrow
T^{\bigtriangledown ^{\tciLaplace }}=0\Leftrightarrow d^{\bigtriangledown
}\tciLaplace =0\Leftrightarrow d^{\bigtriangledown ^{\dag }}g=0$;

$\left( iii\right) $ $d^{\bigtriangledown ^{\dag }}\omega =0\Leftrightarrow
T^{\bigtriangledown }=0\Leftrightarrow d^{\bigtriangledown ^{\tciLaplace
}}\tciLaplace =0\Leftrightarrow d^{\bigtriangledown ^{\ast }}g=0$;

$\left( iv\right) $ $d^{\bigtriangledown ^{\tciLaplace }}\omega
=0\Leftrightarrow T^{\bigtriangledown ^{\ast }}=0\Leftrightarrow
d^{\bigtriangledown ^{\dag }}\tciLaplace =0\Leftrightarrow
d^{\bigtriangledown }g=0$.
\end{proposition}

\begin{proof}
One can write%
\begin{eqnarray*}
&&\left( d^{\bigtriangledown }\omega \right) \left( \xi _{1},\xi _{2},\xi
_{3}\right) \\
&=&\left( \bigtriangledown _{\xi _{1}}\omega \right) \left( \xi _{2},\xi
_{3}\right) -\left( \bigtriangledown _{\xi _{2}}\omega \right) \left( \xi
_{1},\xi _{3}\right) +\omega \left( T^{\bigtriangledown }\left( \xi _{1},\xi
_{2}\right) ,\xi _{3}\right) \\
&=&\xi _{1}\omega \left( \xi _{2},\xi _{3}\right) -\omega \left( \xi
_{2},\bigtriangledown _{\xi _{1}}\xi _{3}\right) -\xi _{2}\omega \left( \xi
_{1},\xi _{3}\right) \\
&&+\omega \left( \xi _{1},\bigtriangledown _{\xi _{2}}\xi _{3}\right)
-\omega \left( \left[ \xi _{1},\xi _{2}\right] ,\xi _{3}\right) \\
&=&\omega \left( \bigtriangledown _{\xi _{1}}^{\dag }\xi _{2},\xi
_{3}\right) -\omega \left( \bigtriangledown _{\xi _{2}}^{\dag }\xi _{1},\xi
_{3}\right) -\omega \left( \left[ \xi _{1},\xi _{2}\right] ,\xi _{3}\right)
\\
&=&\omega \left( T^{\bigtriangledown ^{\dag }}\left( \xi _{1},\xi
_{2}\right) ,\xi _{3}\right) ,
\end{eqnarray*}%
such that $d^{\bigtriangledown }\omega =0\Leftrightarrow T^{\bigtriangledown
^{\dag }}=0$. Moreover, we have%
\begin{eqnarray*}
&&g\left( \left( d^{\bigtriangledown ^{\ast }}\tciLaplace \right) \left( \xi
_{1},\xi _{2}\right) ,\xi _{3}\right) \\
&=&g\left( \left( \bigtriangledown _{\xi _{1}}^{\ast }\tciLaplace \right)
\xi _{2}-\left( \bigtriangledown _{\xi _{2}}^{\ast }\tciLaplace \right) \xi
_{1}+\tciLaplace T^{\bigtriangledown ^{\ast }}\left( \xi _{1},\xi
_{2}\right) ,\xi _{3}\right) \\
&=&g\left( \left( \bigtriangledown _{\xi _{1}}^{\ast }\tciLaplace \right)
\xi _{2},\xi _{3}\right) -g\left( \left( \bigtriangledown _{\xi _{2}}^{\ast
}\tciLaplace \right) \xi _{1},\xi _{3}\right) +g\left( \tciLaplace
T^{\bigtriangledown ^{\ast }}\left( \xi _{1},\xi _{2}\right) ,\xi _{3}\right)
\\
&=&g\left( \bigtriangledown _{\xi _{1}}^{\ast }\tciLaplace \xi _{2},\xi
_{3}\right) -g\left( \bigtriangledown _{\xi _{2}}^{\ast }\tciLaplace \xi
_{1},\xi _{3}\right) -g\left( \tciLaplace \left[ \xi _{1},\xi _{2}\right]
,\xi _{3}\right) \\
&=&\xi _{1}g\left( \tciLaplace \xi _{2},\xi _{3}\right) -g\left( \tciLaplace
\xi _{2},\bigtriangledown _{\xi _{1}}\xi _{3}\right) -\xi _{2}g\left(
\tciLaplace \xi _{1},\xi _{3}\right) \\
&&+g\left( \tciLaplace \xi _{1},\bigtriangledown _{\xi _{2}}\xi _{3}\right)
-g\left( \tciLaplace \left[ \xi _{1},\xi _{2}\right] ,\xi _{3}\right) \\
&=&\xi _{1}\omega \left( \xi _{2},\xi _{3}\right) -\omega \left( \xi
_{2},\bigtriangledown _{\xi _{1}}\xi _{3}\right) -\xi _{2}\omega \left( \xi
_{1},\xi _{3}\right) \\
&&+\omega \left( \xi _{1},\bigtriangledown _{\xi _{2}}\xi _{3}\right)
-\omega \left( \left[ \xi _{1},\xi _{2}\right] ,\xi _{3}\right) \\
&=&\omega \left( T^{\bigtriangledown ^{\dag }}\left( \xi _{1},\xi
_{2}\right) ,\xi _{3}\right) ,
\end{eqnarray*}%
which implies that $T^{\bigtriangledown ^{\dag }}=0\Leftrightarrow
d^{\bigtriangledown ^{\ast }}\tciLaplace =0$. From $(ii)$ of Proposition \ref%
{pro4}, the following last expression is obtained%
\begin{equation*}
d^{\bigtriangledown }\omega =0\Leftrightarrow T^{\bigtriangledown ^{\dag
}}=0\Leftrightarrow d^{\bigtriangledown ^{\ast }}\tciLaplace
=0\Leftrightarrow d^{\bigtriangledown ^{\tciLaplace }}g=0.
\end{equation*}%
With help of Remark \ref{teo1}, the expressions $\left( ii\right) ,\left(
iii\right) $ and $\left( iv\right) $ can be easily proved.
\end{proof}

\begin{remark}
The Proposition \ref{pro5} says to us that the torsion tensor of the $g$%
-conjugate of a linear connection $\nabla $ is always zero on a
quasi-statistical manifold. Also, this proposition gives information about
the integrability of the structure $\tciLaplace $. Suppose that $\left(
\nabla ,\tciLaplace \right) $ is Codazzi-coupled. If the triple $\left( 
\acute{N},g,\nabla ^{\ast }\right) $ is a quasi-statistical manifold, then
the almost complex structure $\tciLaplace $ is integrable.
\end{remark}

\begin{proposition}
Let $\bigtriangledown $ be a linear connection with torsion tensor $%
T^{\nabla }$, $\tciLaplace $ be a almost complex structure and $g$ be a
pseudo-Riemannian metric on $\acute{N}$. If $d^{\bigtriangledown
}\tciLaplace =0$ and $\left( \acute{N},g,\bigtriangledown \right) $ is a
quasi-statistical manifold, then the below equality is satisfied%
\begin{eqnarray*}
&&\left( d^{\bigtriangledown ^{\tciLaplace }}g\right) \left( \xi _{1},\xi
_{2},\xi _{3}\right) +\left( d^{\bigtriangledown ^{\tciLaplace }}g\right)
\left( \xi _{2},\xi _{3},\xi _{1}\right) +\left( d^{\bigtriangledown
^{\tciLaplace }}g\right) \left( \xi _{3},\xi _{1},\xi _{2}\right) \\
&=&g\left( \xi _{2},T^{\bigtriangledown }\left( \xi _{1},\xi _{3}\right)
\right) +g\left( \xi _{3},T^{\bigtriangledown }\left( \xi _{2},\xi
_{1}\right) \right) +g\left( \xi _{1},T^{\bigtriangledown }\left( \xi
_{3},\xi _{2}\right) \right) .
\end{eqnarray*}
\end{proposition}

\begin{proof}
We calculate%
\begin{eqnarray*}
&&\left( d^{\bigtriangledown ^{\tciLaplace }}g\right) \left( \xi _{1},\xi
_{2},\xi _{3}\right) \\
&=&\left( \bigtriangledown _{\xi _{1}}^{\tciLaplace }g\right) \left( \xi
_{2},\xi _{3}\right) -\left( \bigtriangledown _{\xi _{2}}^{\tciLaplace
}g\right) \left( \xi _{1},\xi _{3}\right) +g\left( T^{\bigtriangledown
^{\tciLaplace }}\left( \xi _{1},\xi _{2}\right) ,\xi _{3}\right) \\
&=&\xi _{1}g\left( \xi _{2},\xi _{3}\right) -g\left( \xi
_{2},\bigtriangledown _{\xi _{1}}^{\tciLaplace }\xi _{3}\right) -\xi
_{2}g\left( \xi _{1},\xi _{3}\right) \\
&&+g\left( \xi _{1},\bigtriangledown _{\xi _{2}}^{\tciLaplace }\xi
_{3}\right) -g\left( \left[ \xi _{1},\xi _{2}\right] ,\xi _{3}\right) \\
&=&\xi _{1}g\left( \xi _{2},\xi _{3}\right) -g\left( \xi _{2},\tciLaplace
^{-1}\bigtriangledown _{\xi _{1}}\tciLaplace \xi _{3}\right) -\xi
_{2}g\left( \xi _{1},\xi _{3}\right) \\
&&+g\left( \xi _{1},\tciLaplace ^{-1}\bigtriangledown _{\xi _{2}}\tciLaplace
\xi _{3}\right) -g\left( \left[ \xi _{1},\xi _{2}\right] ,\xi _{3}\right) \\
&=&\xi _{1}g\left( \xi _{2},\xi _{3}\right) -g\left( \xi
_{2},\bigtriangledown _{\xi _{1}}\xi _{3}\right) -g\left( \xi
_{2},\tciLaplace ^{-1}\left( \bigtriangledown _{\xi _{1}}\tciLaplace \right)
\xi _{3}\right) -\xi _{2}g\left( \xi _{1},\xi _{3}\right) \\
&&+g\left( \xi _{1},\bigtriangledown _{\xi _{2}}\xi _{3}\right) +g\left( \xi
_{1},\tciLaplace ^{-1}\left( \bigtriangledown _{\xi _{2}}\tciLaplace \right)
\xi _{3}\right) -g\left( \left[ \xi _{1},\xi _{2}\right] ,\xi _{3}\right) \\
&=&\left( d^{\bigtriangledown }g\right) \left( \xi _{1},\xi _{2},\xi
_{3}\right) -g\left( \xi _{2},\tciLaplace ^{-1}\left( \bigtriangledown _{\xi
_{1}}\tciLaplace \right) \xi _{3}\right) +g\left( \xi _{1},\tciLaplace
^{-1}\left( \bigtriangledown _{\xi _{2}}\tciLaplace \right) \xi _{3}\right) .
\end{eqnarray*}%
The equality $d^{\bigtriangledown }\tciLaplace =0$ implies that 
\begin{eqnarray*}
&&\left( d^{\bigtriangledown ^{\tciLaplace }}g\right) \left( \xi _{1},\xi
_{2},\xi _{3}\right) +\left( d^{\bigtriangledown ^{\tciLaplace }}g\right)
\left( \xi _{2},\xi _{3},\xi _{1}\right) +\left( d^{\bigtriangledown
^{\tciLaplace }}g\right) \left( \xi _{3},\xi _{1},\xi _{2}\right) \\
&=&g\left( \xi _{2},T^{\bigtriangledown }\left( \xi _{1},\xi _{3}\right)
\right) +g\left( \xi _{3},T^{\bigtriangledown }\left( \xi _{2},\xi
_{1}\right) \right) +g\left( \xi _{1},T^{\bigtriangledown }\left( \xi
_{3},\xi _{2}\right) \right) .
\end{eqnarray*}%
Here, we also use $d^{\bigtriangledown }g=0$ and $\tciLaplace ^{2}=-id$.
\end{proof}

\begin{lemma}
\label{lem3}Let $\bigtriangledown $ be a linear connection with torsion
tensor $T^{\nabla }$, $\tciLaplace $ be a almost complex structure and $g$
be a pseudo-Riemannian metric on $\acute{N}$. Let $\left( g,\omega
,\tciLaplace \right) $ be a compatible triple. If $d^{\bigtriangledown
}\tciLaplace =0$ and $\left( \acute{N},g,\bigtriangledown \right) $ is a
quasi-statistical manifold, then $\omega $ is closed, that is, $d\omega =0$.
\end{lemma}

\begin{proof}
From Lemma \ref{lem2}, we obtain%
\begin{eqnarray*}
&&d\omega \left( \xi _{1},\xi _{2},\xi _{3}\right) \\
&=&\left( \bigtriangledown _{\xi _{3}}\omega \right) \left( \xi _{1},\xi
_{2}\right) +\left( \bigtriangledown _{\xi _{1}}\omega \right) \left( \xi
_{2},\xi _{3}\right) +\left( \bigtriangledown _{\xi _{2}}\omega \right)
\left( \xi _{3},\xi _{1}\right) \\
&&+\omega \left( T^{\bigtriangledown }\left( \xi _{1},\xi _{2}\right) ,\xi
_{3}\right) +\omega \left( T^{\bigtriangledown }\left( \xi _{2},\xi
_{3}\right) ,\xi _{1}\right) +\omega \left( T^{\bigtriangledown }\left( \xi
_{3},\xi _{1}\right) ,\xi _{2}\right) \\
&=&-\left( \bigtriangledown _{\xi _{3}}g\right) \left( \xi _{1},\tciLaplace
\xi _{2}\right) -g\left( \xi _{1},\left( \bigtriangledown _{\xi
_{3}}\tciLaplace \right) \xi _{2}\right) -\left( \bigtriangledown _{\xi
_{1}}g\right) \left( \xi _{2},\tciLaplace \xi _{3}\right) \\
&&-g\left( \xi _{2},\left( \bigtriangledown _{\xi _{1}}\tciLaplace \right)
\xi _{3}\right) -\left( \bigtriangledown _{\xi _{2}}g\right) \left( \xi
_{3},\tciLaplace \xi _{1}\right) -g\left( \xi _{3},\left( \bigtriangledown
_{\xi _{2}}\tciLaplace \right) \xi _{1}\right) \\
&&+g\left( \tciLaplace T^{\bigtriangledown }\left( \xi _{1},\xi _{2}\right)
,\xi _{3}\right) +g\left( \tciLaplace T^{\bigtriangledown }\left( \xi
_{2},\xi _{3}\right) ,\xi _{1}\right) +g\left( \tciLaplace
T^{\bigtriangledown }\left( \xi _{3},\xi _{1}\right) ,\xi _{2}\right) \\
&=&-\left( \bigtriangledown _{\xi _{1}}g\right) \left( \xi _{3},\tciLaplace
\xi _{2}\right) -\left( \bigtriangledown _{\xi _{2}}g\right) \left( \xi
_{1},\tciLaplace \xi _{3}\right) -\left( \bigtriangledown _{\xi
_{3}}g\right) \left( \xi _{2},\tciLaplace \xi _{1}\right) \\
&&-g\left( \xi _{1},\left( \bigtriangledown _{\xi _{3}}\tciLaplace \right)
\xi _{2}\right) -g\left( \xi _{2},\left( \bigtriangledown _{\xi
_{1}}\tciLaplace \right) \xi _{3}\right) -g\left( \xi _{3},\left(
\bigtriangledown _{\xi _{2}}\tciLaplace \right) \xi _{1}\right)
\end{eqnarray*}%
and 
\begin{eqnarray*}
d\omega \left( \xi _{3},\xi _{2},\xi _{1}\right) &=&-\left( \bigtriangledown
_{\xi _{3}}g\right) \left( \xi _{1},\tciLaplace \xi _{2}\right) -\left(
\bigtriangledown _{\xi _{2}}g\right) \left( \xi _{3},\tciLaplace \xi
_{1}\right) -\left( \bigtriangledown _{\xi _{1}}g\right) \left( \xi
_{2},\tciLaplace \xi _{3}\right) \\
&&-g\left( \xi _{3},\left( \bigtriangledown _{\xi _{1}}\tciLaplace \right)
\xi _{2}\right) -g\left( \xi _{2},\left( \bigtriangledown _{\xi
_{3}}\tciLaplace \right) \xi _{1}\right) -g\left( \xi _{1},\left(
\bigtriangledown _{\xi _{2}}\tciLaplace \right) \xi _{3}\right) .
\end{eqnarray*}%
Thus, we have the following%
\begin{eqnarray*}
&&d\omega \left( \xi _{1},\xi _{2},\xi _{3}\right) -d\omega \left( \xi
_{3},\xi _{2},\xi _{1}\right) \\
&=&-\left( \bigtriangledown _{\xi _{1}}g\right) \left( \xi _{3},\tciLaplace
\xi _{2}\right) -\left( \bigtriangledown _{\xi _{2}}g\right) \left( \xi
_{1},\tciLaplace \xi _{3}\right) -\left( \bigtriangledown _{\xi
_{3}}g\right) \left( \xi _{2},\tciLaplace \xi _{1}\right) \\
&&-g\left( \xi _{1},\left( \bigtriangledown _{\xi _{3}}\tciLaplace \right)
\xi _{2}\right) -g\left( \xi _{2},\left( \bigtriangledown _{\xi
_{1}}\tciLaplace \right) \xi _{3}\right) -g\left( \xi _{3},\left(
\bigtriangledown _{\xi _{2}}\tciLaplace \right) \xi _{1}\right) \\
&&+\left( \bigtriangledown _{\xi _{3}}g\right) \left( \xi _{1},\tciLaplace
\xi _{2}\right) +\left( \bigtriangledown _{\xi _{2}}g\right) \left( \xi
_{3},\tciLaplace \xi _{1}\right) +\left( \bigtriangledown _{\xi
_{1}}g\right) \left( \xi _{2},\tciLaplace \xi _{3}\right) \\
&&+g\left( \xi _{3},\left( \bigtriangledown _{\xi _{1}}\tciLaplace \right)
\xi _{2}\right) +g\left( \xi _{2},\left( \bigtriangledown _{\xi
_{3}}\tciLaplace \right) \xi _{1}\right) +g\left( \xi _{1},\left(
\bigtriangledown _{\xi _{2}}\tciLaplace \right) \xi _{3}\right) \\
&=&-\left( \bigtriangledown _{\xi _{1}}g\right) \left( \xi _{3},\tciLaplace
\xi _{2}\right) +\left( \bigtriangledown _{\xi _{3}}g\right) \left( \xi
_{1},\tciLaplace \xi _{2}\right) \\
&&+g\left( \xi _{2},\left( \bigtriangledown _{\xi _{3}}\tciLaplace \right)
\xi _{1}-\left( \bigtriangledown _{\xi _{1}}\tciLaplace \right) \xi
_{3}\right) \\
&&-\left( \bigtriangledown _{\xi _{2}}g\right) \left( \xi _{1},\tciLaplace
\xi _{3}\right) +\left( \bigtriangledown _{\xi _{1}}g\right) \left( \xi
_{2},\tciLaplace \xi _{3}\right) \\
&&+g\left( \xi _{3},\left( \bigtriangledown _{\xi _{1}}\tciLaplace \right)
\xi _{2}-\left( \bigtriangledown _{\xi _{2}}\tciLaplace \right) \xi
_{1}\right) \\
&&-\left( \bigtriangledown _{\xi _{3}}g\right) \left( \xi _{2},\tciLaplace
\xi _{1}\right) +\left( \bigtriangledown _{\xi _{2}}g\right) \left( \xi
_{3},\tciLaplace \xi _{1}\right) \\
&&+g\left( \xi _{1},\left( \bigtriangledown _{\xi _{2}}\tciLaplace \right)
\xi _{3}-\left( \bigtriangledown _{\xi _{3}}\tciLaplace \right) \xi
_{2}\right) \\
&=&-\left( \bigtriangledown _{\xi _{1}}g\right) \left( \xi _{3},\tciLaplace
\xi _{2}\right) +\left( \bigtriangledown _{\xi _{3}}g\right) \left( \xi
_{1},\tciLaplace \xi _{2}\right) +g\left( \xi _{2},\tciLaplace
T^{\bigtriangledown }\left( \xi _{1},\xi _{3}\right) \right) \\
&&-\left( \bigtriangledown _{\xi _{2}}g\right) \left( \xi _{1},\tciLaplace
\xi _{3}\right) +\left( \bigtriangledown _{\xi _{1}}g\right) \left( \xi
_{2},\tciLaplace \xi _{3}\right) +g\left( \xi _{3},\tciLaplace
T^{\bigtriangledown }\left( \xi _{2},\xi _{1}\right) \right) \\
&&-\left( \bigtriangledown _{\xi _{3}}g\right) \left( \xi _{2},\tciLaplace
\xi _{1}\right) +\left( \bigtriangledown _{\xi _{2}}g\right) \left( \xi
_{3},\tciLaplace \xi _{1}\right) +g\left( \xi _{1},\tciLaplace
T^{\bigtriangledown }\left( \xi _{3},\xi _{2}\right) \right) \\
&=&-\left( d^{\bigtriangledown }g\right) \left( \xi _{1},\xi
_{3},\tciLaplace \xi _{2}\right) -\left( d^{\bigtriangledown }g\right)
\left( \xi _{2},\xi _{1},\tciLaplace \xi _{3}\right) -\left(
d^{\bigtriangledown }g\right) \left( \xi _{3},\xi _{2},\tciLaplace \xi
_{1}\right) \\
&=&0,
\end{eqnarray*}%
such that $d\omega \left( \xi _{1},\xi _{2},\xi _{3}\right) =d\omega \left(
\xi _{3},\xi _{2},\xi _{1}\right) $. Since $d\omega $ is totally
skew-symmetric, we conclude $d\omega =0$, i.e., $\omega $ is closed.
\end{proof}

We are now ready to introduce our first main theorem. As a corollary of
Proposition \ref{pro2} and Lemma \ref{lem3}, we have the following.

\begin{theorem}
\label{teo2}Let $\bigtriangledown $ be a linear connection with torsion
tensor $T^{\nabla }$, $g$ be a pseudo-Riemannian metric and $\tciLaplace $
be an almost complex structure on $\acute{N}$, and $\left( g,\omega
,\tciLaplace \right) $ be a compatible triple. Assume that the torsion
tensor $T^{\bigtriangledown }$ of $\nabla $ satisfies the
torsion-compatibility condition. If $d^{\bigtriangledown }\tciLaplace =0$
and $\left( \acute{N},g,\bigtriangledown \right) $ is a quasi-statistical
manifold, then $\left( \acute{N},g,\nabla ,\tciLaplace \right) $ is a K\"{a}%
hler manifold.
\end{theorem}

\begin{remark}
We know that for any statistical manifold $\left( \acute{N},g,\nabla \right) 
$, if there exists Codazzi couplings of $\nabla $ with an almost complex
structure $\tciLaplace $, $\left( \acute{N},g,\nabla ,\tciLaplace \right) $
is a K\"{a}hler manifold (see Theorem 3.2 in \cite{Fei} ). Theorem \ref{teo2}
says that an alternative classification can be made for K\"{a}hler manifolds
by taking a quasi-statistical structure instead of a statistical structure.
That is, to make such a classification, the torsion tensor of a linear
connection $\nabla $ need not be zero.
\end{remark}

\begin{theorem}
Let $\bigtriangledown $ be a linear connection with torsion tensor $%
T^{\nabla }$ on $\acute{N}$, and $\left( g,\omega ,\tciLaplace \right) $ be
a compatible triple. Then, for the following three statements regarding any
compatible triple $\left( g,\omega ,\tciLaplace \right) $, any two imply the
third

$\left( i\right) $ $\left( \acute{N},\bigtriangledown ,g\right) $ is a
quasi-statistical manifold;

$\left( ii\right) $ $d^{\bigtriangledown }\tciLaplace =0$, that is, $%
\tciLaplace $ is $d^{\nabla }$-closed;

$\left( iii\right) $ $\bigtriangledown ^{\ast }\omega =0$.
\end{theorem}

\begin{proof}
Assume that $d^{\bigtriangledown }g=0$ and $d^{\bigtriangledown }\tciLaplace
=0$. From Lemma \ref{lem3}, we have%
\begin{eqnarray*}
&&d\omega \left( \xi _{1},\xi _{2},\xi _{3}\right) \\
&=&-\left( \bigtriangledown _{\xi _{3}}g\right) \left( \xi _{1},\tciLaplace
\xi _{2}\right) -g\left( \xi _{1},\left( \bigtriangledown _{\xi
_{3}}\tciLaplace \right) \xi _{2}\right) -\left( \bigtriangledown _{\xi
_{1}}g\right) \left( \xi _{2},\tciLaplace \xi _{3}\right) \\
&&-g\left( \xi _{2},\left( \bigtriangledown _{\xi _{1}}\tciLaplace \right)
\xi _{3}\right) -\left( \bigtriangledown _{\xi _{2}}g\right) \left( \xi
_{3},\tciLaplace \xi _{1}\right) -g\left( \xi _{3},\left( \bigtriangledown
_{\xi _{2}}\tciLaplace \right) \xi _{1}\right) \\
&&+g\left( \tciLaplace T^{\bigtriangledown }\left( \xi _{1},\xi _{2}\right)
,\xi _{3}\right) +g\left( \tciLaplace T^{\bigtriangledown }\left( \xi
_{2},\xi _{3}\right) ,\xi _{1}\right) +g\left( \tciLaplace
T^{\bigtriangledown }\left( \xi _{3},\xi _{1}\right) ,\xi _{2}\right) \\
&=&-\left( \bigtriangledown _{\xi _{3}}g\right) \left( \xi _{1},\tciLaplace
\xi _{2}\right) -g\left( \xi _{1},\left( \bigtriangledown _{\xi
_{3}}\tciLaplace \right) \xi _{2}\right) -\left( \bigtriangledown _{\xi
_{1}}g\right) \left( \xi _{2},\tciLaplace \xi _{3}\right) \\
&&-g\left( \xi _{2},\left( \bigtriangledown _{\xi _{1}}\tciLaplace \right)
\xi _{3}\right) -\left( \bigtriangledown _{\xi _{2}}g\right) \left( \xi
_{3},\tciLaplace \xi _{1}\right) -g\left( \xi _{3},\left( \bigtriangledown
_{\xi _{2}}\tciLaplace \right) \xi _{1}\right) \\
&&-g\left( T^{\bigtriangledown }\left( \xi _{1},\xi _{2}\right) ,\tciLaplace
\xi _{3}\right) -g\left( T^{\bigtriangledown }\left( \xi _{2},\xi
_{3}\right) ,\tciLaplace \xi _{1}\right) -g\left( T^{\bigtriangledown
}\left( \xi _{3},\xi _{1}\right) ,\tciLaplace \xi _{2}\right) \\
&=&-\left( \bigtriangledown _{\xi _{1}}g\right) \left( \xi _{3},\tciLaplace
\xi _{2}\right) -\left( \bigtriangledown _{\xi _{2}}g\right) \left( \xi
_{1},\tciLaplace \xi _{3}\right) -\left( \bigtriangledown _{\xi
_{3}}g\right) \left( \xi _{2},\tciLaplace \xi _{1}\right) \\
&&-g\left( \xi _{1},\left( \bigtriangledown _{\xi _{3}}\tciLaplace \right)
\xi _{2}\right) -g\left( \xi _{2},\left( \bigtriangledown _{\xi
_{1}}\tciLaplace \right) \xi _{3}\right) -g\left( \xi _{3},\left(
\bigtriangledown _{\xi _{2}}\tciLaplace \right) \xi _{1}\right) \\
&=&-\left( \bigtriangledown _{\xi _{1}}g\right) \left( \xi _{3},\tciLaplace
\xi _{2}\right) -\left( \bigtriangledown _{\xi _{2}}g\right) \left( \xi
_{1},\tciLaplace \xi _{3}\right) -\left( \bigtriangledown _{\xi
_{3}}g\right) \left( \xi _{2},\tciLaplace \xi _{1}\right) \\
&&-g\left( \xi _{1},\tciLaplace T^{\bigtriangledown }\left( \xi _{2},\xi
_{3}\right) +\left( \bigtriangledown _{\xi _{2}}\tciLaplace \right) \xi
_{3}\right) \\
&&-g\left( \xi _{2},\tciLaplace T^{\bigtriangledown }\left( \xi _{3},\xi
_{1}\right) +\left( \bigtriangledown _{\xi _{3}}\tciLaplace \right) \xi
_{1}\right) \\
&&-g\left( \xi _{3},\tciLaplace T^{\bigtriangledown }\left( \xi _{1},\xi
_{2}\right) +\left( \bigtriangledown _{\xi _{1}}\tciLaplace \right) \xi
_{2}\right) \\
&=&-\left( \bigtriangledown _{\xi _{3}}g\right) \left( \xi _{1},\tciLaplace
\xi _{2}\right) -\left( \bigtriangledown _{\xi _{1}}g\right) \left( \xi
_{2},\tciLaplace \xi _{3}\right) -\left( \bigtriangledown _{\xi
_{2}}g\right) \left( \xi _{3},\tciLaplace \xi _{1}\right) \\
&&-g\left( \xi _{1},\left( \bigtriangledown _{\xi _{2}}\tciLaplace \right)
\xi _{3}\right) -g\left( \xi _{2},\left( \bigtriangledown _{\xi
_{3}}\tciLaplace \right) \xi _{1}\right) -g\left( \xi _{3},\left(
\bigtriangledown _{\xi _{1}}\tciLaplace \right) \xi _{2}\right) .
\end{eqnarray*}%
Besides, due to skew-symmetric of $\omega $, we can write 
\begin{eqnarray*}
&&\left( \bigtriangledown _{\xi _{3}}\omega \right) \left( \xi _{1},\xi
_{2}\right) +\left( \bigtriangledown _{\xi _{3}}\omega \right) \left( \xi
_{2},\xi _{1}\right) \\
&=&-\left( \bigtriangledown _{\xi _{3}}g\right) \left( \xi _{1},\tciLaplace
\xi _{2}\right) -g\left( \xi _{1},\left( \bigtriangledown _{\xi
_{3}}\tciLaplace \right) \xi _{2}\right) \\
&&-\left( \bigtriangledown _{\xi _{3}}g\right) \left( \xi _{2},\tciLaplace
\xi _{1}\right) -g\left( \xi _{2},\left( \bigtriangledown _{\xi
_{3}}\tciLaplace \right) \xi _{1}\right) \\
&=&-\left( \bigtriangledown _{\xi _{3}}g\right) \left( \xi _{1},\tciLaplace
\xi _{2}\right) -g\left( \xi _{1},\left( \bigtriangledown _{\xi
_{3}}\tciLaplace \right) \xi _{2}\right) -\left( \bigtriangledown _{\xi
_{2}}g\right) \left( \xi _{3},\tciLaplace \xi _{1}\right) \\
&&-g\left( T^{\bigtriangledown }\left( \xi _{2},\xi _{3}\right) ,\tciLaplace
\xi _{1}\right) -g\left( \xi _{2},\left( \bigtriangledown _{\xi
_{3}}\tciLaplace \right) \xi _{1}\right) \\
&=&-\left( \bigtriangledown _{\xi _{3}}g\right) \left( \xi _{1},\tciLaplace
\xi _{2}\right) -\left( \bigtriangledown _{\xi _{2}}g\right) \left( \xi
_{3},\tciLaplace \xi _{1}\right) \\
&&-g\left( \xi _{2},\left( \bigtriangledown _{\xi _{3}}\tciLaplace \right)
\xi _{1}\right) -g\left( \xi _{1},\left( \bigtriangledown _{\xi
_{2}}\tciLaplace \right) \xi _{3}\right) \\
&=&0.
\end{eqnarray*}%
With these relations, we get%
\begin{equation*}
d\omega \left( \xi _{1},\xi _{2},\xi _{3}\right) =-\left( \bigtriangledown
_{\xi _{1}}g\right) \left( \xi _{2},\tciLaplace \xi _{3}\right) -g\left( \xi
_{3},\left( \bigtriangledown _{\xi _{1}}\tciLaplace \right) \xi _{2}\right)
=-\left( \bigtriangledown _{\xi _{1}}^{\ast }\omega \right) \left( \xi
_{2},\xi _{3}\right) =0,
\end{equation*}%
that is, $\bigtriangledown ^{\ast }\omega =0$.\newline

Next, let us suppose that $d^{\bigtriangledown }g=0$ and $\bigtriangledown
^{\ast }\omega =0$. Thus, we get%
\begin{eqnarray*}
\left( d^{\bigtriangledown }g\right) \left( \xi _{1},\xi _{2},\tciLaplace
\xi _{3}\right) &=&\left( \bigtriangledown _{\xi _{1}}g\right) \left( \xi
_{2},\tciLaplace \xi _{3}\right) -\left( \bigtriangledown _{\xi
_{2}}g\right) \left( \xi _{1},\tciLaplace \xi _{3}\right) \\
&&+g\left( T^{\bigtriangledown }\left( \xi _{1},\xi _{2}\right) ,\tciLaplace
\xi _{3}\right) \\
&=&0
\end{eqnarray*}%
and%
\begin{eqnarray*}
\left( \bigtriangledown _{\xi _{1}}^{\ast }\omega \right) \left( \xi
_{2},\xi _{3}\right) &=&\left( \bigtriangledown _{\xi _{1}}g\right) \left(
\xi _{2},\tciLaplace \xi _{3}\right) +g\left( \xi _{3},\left(
\bigtriangledown _{\xi _{1}}\tciLaplace \right) \xi _{2}\right) \\
&=&0
\end{eqnarray*}%
from which we immediately see that%
\begin{eqnarray*}
&&\left( d^{\bigtriangledown }g\right) \left( \xi _{1},\xi _{2},\tciLaplace
\xi _{3}\right) \\
&=&-g\left( \xi _{3},\left( \bigtriangledown _{\xi _{1}}\tciLaplace \right)
\xi _{2}\right) +g\left( \xi _{3},\left( \bigtriangledown _{\xi
_{2}}\tciLaplace \right) \xi _{1}\right) -g\left( \tciLaplace
T^{\bigtriangledown }\left( \xi _{1},\xi _{2}\right) ,\xi _{3}\right) \\
&=&-g\left( \xi _{3},\left( d^{\bigtriangledown }\tciLaplace \right) \left(
\xi _{1},\xi _{2}\right) \right) =0,
\end{eqnarray*}%
such that $d^{\bigtriangledown }\tciLaplace =0$. It is easy to see that if $%
d^{\bigtriangledown }\tciLaplace =0$ and $\bigtriangledown ^{\ast }\omega =0$%
, then $d^{\bigtriangledown }g=0$.
\end{proof}

Now we turn our attention to the linear connection $\widetilde{%
\bigtriangledown }$ which is the average of a linear connection and its $%
\tciLaplace $-conjugate connection such that $\widetilde{\bigtriangledown }=%
\frac{1}{2}\left( \bigtriangledown ^{\tciLaplace }+\bigtriangledown \right) $%
. The connection $\widetilde{\bigtriangledown }$ is a complex connection,
that is, $\widetilde{\bigtriangledown }\tciLaplace =0$ \cite{Grigorian}.

\begin{proposition}
\label{GAD15} Let $\left( \acute{N},g\right) $ be a pseudo-Riemannian
manifold, $\bigtriangledown $ be an arbitrary linear connection with torsion
tensor $T^{\nabla }$, $\bigtriangledown ^{\ast }$be the $g$-conjugate
connection of $\bigtriangledown $ and $\tciLaplace $ be an almost complex
structure that is compatible with $g$. Assume that $\left( \bigtriangledown
^{\ast },\tciLaplace \right) $ is Codazzi-coupled. $\left( \widetilde{%
\bigtriangledown },g\right) $ is a quasi-statistical manifold if and only if 
$\left( \bigtriangledown ,g\right) $ is a quasi-statistical manifold, where $%
\widetilde{\bigtriangledown }=\frac{1}{2}\left( \bigtriangledown
^{\tciLaplace }+\bigtriangledown \right) $ and $\bigtriangledown
^{\tciLaplace }$ is the $\tciLaplace $-conjugate connection of $\nabla $.
\end{proposition}

\begin{proof}
We have%
\begin{eqnarray*}
&&\left( d^{\widetilde{\bigtriangledown }}g\right) \left( \xi _{1},\xi
_{2},\xi _{3}\right) \\
&=&\left( \widetilde{\bigtriangledown }_{\xi _{1}}g\right) \left( \xi
_{2},\xi _{3}\right) -\left( \widetilde{\bigtriangledown }_{\xi
_{2}}g\right) \left( \xi _{1},\xi _{3}\right) +g\left( T^{\widetilde{%
\bigtriangledown }}\left( \xi _{1},\xi _{2}\right) ,\xi _{3}\right) \\
&=&\frac{1}{2}\left( \bigtriangledown _{\xi _{1}}g\right) \left( \xi
_{2},\xi _{3}\right) +\frac{1}{2}\left( \bigtriangledown _{\xi
_{1}}^{\tciLaplace }g\right) \left( \xi _{2},\xi _{3}\right) -\frac{1}{2}%
\left( \bigtriangledown _{\xi _{2}}g\right) \left( \xi _{1},\xi _{3}\right)
\\
&&-\frac{1}{2}\left( \bigtriangledown _{\xi _{2}}^{\tciLaplace }g\right)
\left( \xi _{1},\xi _{3}\right) +\frac{1}{2}g\left( T^{\bigtriangledown
}\left( \xi _{1},\xi _{2}\right) ,\xi _{3}\right) \\
&&+\frac{1}{2}g\left( T^{\bigtriangledown ^{\tciLaplace }}\left( \xi
_{1},\xi _{2}\right) ,\xi _{3}\right) .
\end{eqnarray*}%
On considering the following equalities 
\begin{eqnarray*}
\left( \bigtriangledown _{\xi _{1}}^{\tciLaplace }g\right) \left( \xi
_{2},\xi _{3}\right) &=&\left( \bigtriangledown _{\xi _{1}}g\right) \left(
\xi _{2},\xi _{3}\right) -g\left( \tciLaplace ^{-1}\left( \bigtriangledown
_{\xi _{1}}\tciLaplace \right) \xi _{2},\xi _{3}\right) \\
&&-g\left( \xi _{2},\tciLaplace ^{-1}\left( \bigtriangledown _{\xi
_{1}}\tciLaplace \right) \xi _{3}\right)
\end{eqnarray*}%
and 
\begin{eqnarray*}
g\left( T^{\bigtriangledown ^{\tciLaplace }}\left( \xi _{1},\xi _{2}\right)
,\xi _{3}\right) &=&g\left( T^{\bigtriangledown }\left( \xi _{1},\xi
_{2}\right) ,\xi _{3}\right) +g\left( \tciLaplace ^{-1}\left(
\bigtriangledown _{\xi _{1}}\tciLaplace \right) \xi _{2},\xi _{3}\right) \\
&&-g\left( \tciLaplace ^{-1}\left( \bigtriangledown _{\xi _{2}}\tciLaplace
\right) \xi _{1},\xi _{3}\right) ,
\end{eqnarray*}%
we obtain%
\begin{eqnarray*}
&&\left( d^{\widetilde{\bigtriangledown }}g\right) \left( \xi _{1},\xi
_{2},\xi _{3}\right) \\
&=&\left( d^{\bigtriangledown }g\right) \left( \xi _{1},\xi _{2},\xi
_{3}\right) -\frac{1}{2}g\left( \xi _{2},\tciLaplace ^{-1}\left(
\bigtriangledown _{\xi _{1}}\tciLaplace \right) \xi _{3}\right) \\
&&+\frac{1}{2}g\left( \xi _{1},\tciLaplace ^{-1}\left( \bigtriangledown
_{\xi _{2}}\tciLaplace \right) \xi _{3}\right) .
\end{eqnarray*}%
To complete the proof, we need to note $g\left( \xi _{1},\tciLaplace
^{-1}\left( \bigtriangledown _{\xi _{2}}\tciLaplace \right) \xi _{3}\right)
=-g\left( \left( \bigtriangledown _{\xi _{2}}^{\ast }\tciLaplace \right) \xi
_{1},\tciLaplace \xi _{3}\right) $. From the hypothesis, we get $\left( d^{%
\widetilde{\bigtriangledown }}g\right) \left( \xi _{1},\xi _{2},\xi
_{3}\right) =\left( d^{\bigtriangledown }g\right) \left( \xi _{1},\xi
_{2},\xi _{3}\right) $.
\end{proof}

Via $(vi)$ of Proposition \ref{pro3}, $(i)$ and $(iv)$ of Proposition \ref%
{pro5} and Proposition \ref{GAD15}, we obtain the following result.

\begin{corollary}
Assume that $\left( \nabla ^{\ast },\tciLaplace \right) $ is
Codazzi-coupled. Then, $\left( \widetilde{\bigtriangledown },g\right) $ is a
quasi-statistical structure if and only if the torsion tensors of $\nabla
^{\ast }$ and $\nabla ^{\dag }$ are zero.
\end{corollary}

\begin{proposition}
\label{GAD16}Let $\left( \acute{N},g\right) $ be a pseudo-Riemannian
manifold equipped with a pseudo-Riemannian metric $g$ and a linear
connection $\bigtriangledown $ with torsion tensor $T^{\nabla }$. Let $%
\left( g,\omega ,\tciLaplace \right) $ be a compatible triple, and $%
\bigtriangledown ^{\ast },\bigtriangledown ^{\dag }$ denote, respectively, $%
g $-conjugation, $\omega $-conjugation of the linear connection $%
\bigtriangledown $. Then, $\left( \widetilde{\bigtriangledown },g\right) $
is a quasi-statistical manifold if and only if $T^{\bigtriangledown ^{\ast
}}=-T^{\bigtriangledown ^{\dag }}$, where $\widetilde{\bigtriangledown }=%
\frac{1}{2}\left( \bigtriangledown ^{\tciLaplace }+\bigtriangledown \right) $%
.
\end{proposition}

\begin{proof}
Considering the definition of $d^{\widetilde{\bigtriangledown }}g$, we have%
\begin{eqnarray*}
&&\left( d^{\widetilde{\bigtriangledown }}g\right) \left( \xi _{1},\xi
_{2},\xi _{3}\right) \\
&=&\left( \widetilde{\bigtriangledown }_{\xi _{1}}g\right) \left( \xi
_{2},\xi _{3}\right) -\left( \widetilde{\bigtriangledown }_{\xi
_{2}}g\right) \left( \xi _{1},\xi _{3}\right) +g\left( T^{\widetilde{%
\bigtriangledown }}\left( \xi _{1},\xi _{2}\right) ,\xi _{3}\right) \\
&=&\xi _{1}g\left( \xi _{2},\xi _{3}\right) -g\left( \xi _{2},\widetilde{%
\bigtriangledown }_{\xi _{1}}\xi _{3}\right) -\xi _{2}g\left( \xi _{1},\xi
_{3}\right) \\
&&+g\left( \xi _{1},\widetilde{\bigtriangledown }_{\xi _{2}}\xi _{3}\right)
-g\left( \left[ \xi _{1},\xi _{2}\right] ,\xi _{3}\right) \\
&=&\frac{1}{2}g\left( \bigtriangledown _{\xi _{1}}^{\ast }\xi _{2},\xi
_{3}\right) +\frac{1}{2}g\left( \bigtriangledown _{\xi _{1}}^{\dag }\xi
_{2},\xi _{3}\right) -\frac{1}{2}g\left( \bigtriangledown _{\xi _{2}}^{\ast
}\xi _{1},\xi _{3}\right) \\
&&-\frac{1}{2}g\left( \bigtriangledown _{\xi _{2}}^{\dag }\xi _{1},\xi
_{3}\right) -g\left( \left[ \xi _{1},\xi _{2}\right] ,\xi _{3}\right) \\
&=&\frac{1}{2}g\left( T^{\bigtriangledown ^{\ast }}\left( \xi _{1},\xi
_{2}\right) +T^{\bigtriangledown ^{\dag }}\left( \xi _{1},\xi _{2}\right)
,\xi _{3}\right) .
\end{eqnarray*}
\end{proof}

\begin{proposition}
\label{GAD17} Let $\left( \acute{N},g\right) $ be a pseudo-Riemannian
manifold equipped with a pseudo-Riemannian metric $g$ and a linear
connection $\bigtriangledown $ with torsion tensor $T^{\nabla }$. Let $%
\left( g,\omega ,\tciLaplace \right) $ be a compatible triple, and $%
\bigtriangledown ^{\ast },\bigtriangledown ^{\dag }$ denote, respectively, $%
g $-conjugation, $\omega $-conjugation of an arbitrary linear connection $%
\bigtriangledown $. Then, $T^{\bigtriangledown ^{\ast
}}=-T^{\bigtriangledown ^{\dag }}$ if and only if $d^{\bigtriangledown
^{\ast }}\tciLaplace =-d^{\bigtriangledown ^{\dag }}\tciLaplace $.
\end{proposition}

\begin{proof}
We calculate%
\begin{eqnarray*}
&&g\left( \left( d^{\bigtriangledown ^{\ast }}\tciLaplace \right) \left( \xi
_{1},\xi _{2}\right) +\left( d^{\bigtriangledown ^{\dag }}\tciLaplace
\right) \left( \xi _{1},\xi _{2}\right) ,\xi _{3}\right) \\
&=&g\left( \left( d^{\bigtriangledown ^{\ast }}\tciLaplace \right) \left(
\xi _{1},\xi _{2}\right) ,\xi _{3}\right) +g\left( \left(
d^{\bigtriangledown ^{\dag }}\tciLaplace \right) \left( \xi _{1},\xi
_{2}\right) ,\xi _{3}\right) \\
&=&g\left( \left( \bigtriangledown _{\xi _{1}}^{\ast }\tciLaplace \right)
\xi _{2}-\left( \bigtriangledown _{\xi _{2}}^{\ast }\tciLaplace \right) \xi
_{1}+\tciLaplace T^{\bigtriangledown ^{\ast }}\left( \xi _{1},\xi
_{2}\right) ,\xi _{3}\right) \\
&&+g\left( \left( \bigtriangledown _{\xi _{1}}^{\dag }\tciLaplace \right)
\xi _{2}-\left( \bigtriangledown _{\xi _{2}}^{\dag }\tciLaplace \right) \xi
_{1}+\tciLaplace T^{\bigtriangledown ^{\dag }}\left( \xi _{1},\xi
_{2}\right) ,\xi _{3}\right) \\
&=&g\left( \bigtriangledown _{\xi _{1}}^{\ast }\tciLaplace \xi _{2},\xi
_{3}\right) -g\left( \bigtriangledown _{\xi _{2}}^{\ast }\tciLaplace \xi
_{1},\xi _{3}\right) -g\left( \tciLaplace \left[ \xi _{1},\xi _{2}\right]
,\xi _{3}\right) \\
&&+g\left( \bigtriangledown _{\xi _{1}}^{\dag }\tciLaplace \xi _{2},\xi
_{3}\right) -g\left( \bigtriangledown _{\xi _{2}}^{\dag }\tciLaplace \xi
_{1},\xi _{3}\right) -g\left( \tciLaplace \left[ \xi _{1},\xi _{2}\right]
,\xi _{3}\right) \\
&=&\omega \left( T^{\bigtriangledown ^{\dag }}\left( \xi _{1},\xi
_{2}\right) +T^{\bigtriangledown ^{\ast }}\left( \xi _{1},\xi _{2}\right)
,\xi _{3}\right) .
\end{eqnarray*}
\end{proof}

Propositions \ref{GAD16} and \ref{GAD17} immediately give the following.

\begin{corollary}
$\left( \widetilde{\bigtriangledown },g\right) $ is a quasi-statistical
structure $\Leftrightarrow T^{\bigtriangledown ^{\ast
}}=-T^{\bigtriangledown ^{\dag }}\Leftrightarrow d^{\bigtriangledown ^{\ast
}}\tciLaplace =-d^{\bigtriangledown ^{\dag }}\tciLaplace $.
\end{corollary}

\section{Quasi Statistical Structures with the anti-Hermitian metric $h$}

In this section, we will investigate the properties of quasi-statistical
manifolds by taking anti-Hermitian metric $h$ instead of the Hermitian
metric $g$. Moreover, considering any linear connection $\nabla $ with
torsion tensor $T^{\nabla }$ instead of the Levi-Civita connection $\nabla
^{h}$ of the anti-Hermitian metric $h$, we will show that the anti-K\"{a}%
hler and quasi-K\"{a}hler-Norden manifolds can be classified under certain
conditions.

\begin{definition}
On given a pseudo-Riemannian manifold $\left( \acute{N},h\right) $ endowed
with an almost complex structure $\tciLaplace $, then the triple $\left( 
\acute{N},\tciLaplace ,h\right) $ is called an almost anti-Hermitian
manifold (or Norden manifold) if%
\begin{equation*}
h\left( \tciLaplace \xi _{1},\xi _{2}\right) =h\left( \xi _{1},\tciLaplace
\xi _{2}\right)
\end{equation*}%
for any vector fields $\xi _{1}$ and $\xi _{2}$ on $\acute{N}$, where the
signature of $h$ is $\left( n,n\right) $, that is, $h$ is a neutral metric.
If the structure $\tciLaplace $ is integrable, then the triple $\left( 
\acute{N},\tciLaplace ,h\right) $ is called an anti-Hermitian manifold or
complex Norden manifold. Also, the twin anti-Hermitian metric is defined by%
\begin{equation*}
\hslash \left( \xi _{1},\xi _{2}\right) =h\left( \tciLaplace \xi _{1},\xi
_{2}\right)
\end{equation*}%
for any vector fields $\xi _{1}$ and $\xi _{2}$. An anti-K\"{a}hler manifold
is an almost anti-Hermitian manifold such that $\nabla ^{h}\tciLaplace =0$,
where $\nabla ^{h}$ is the Levi-Civita connection of the pseudo-Riemannian
manifold $\left( \acute{N},h\right) $ \cite{Ganchev,Iscan,Manev}.
\end{definition}

\bigskip

The Tachibana operator on an almost anti-Hermitian manifold $\left( \acute{N}%
,\tciLaplace ,h\right) $ 
\begin{equation*}
\Phi _{\tciLaplace }:\Im _{2}^{0}(\acute{N})\longrightarrow \Im _{0}^{0}(%
\acute{N})
\end{equation*}%
which is defined from the set of all $(0,2)$-tensor fields ($\Im _{2}^{0}(%
\acute{N})$) into the set of all $(0,3)$-tensor fields ($\Im _{2}^{0}(\acute{%
N})$) on $\acute{N}$ is given by \cite{SalimovK,Tachibana} 
\begin{eqnarray*}
\left( \Phi _{\tciLaplace }h\right) \left( \xi _{1},\xi _{2},\xi _{3}\right)
&=&\tciLaplace \xi _{1}h\left( \xi _{2},\xi _{3}\right) -\xi _{1}h\left(
\tciLaplace \xi _{2},\xi _{3}\right) \\
&&+h\left( \left( L_{\xi _{2}}\tciLaplace \right) \xi _{1},\xi _{3}\right)
+h\left( \xi _{2},\left( L_{\xi _{3}}\tciLaplace \right) \xi _{1}\right) ,
\end{eqnarray*}%
where $\left( L_{\xi _{1}}\tciLaplace \right) \xi _{2}=\left[ \xi
_{1},\tciLaplace \xi _{2}\right] -\tciLaplace \left[ \xi _{1},\xi _{2}\right]
$.

\begin{definition}
\cite{Manev} An almost anti-Hermitian manifold is called a quasi-K\"{a}%
hler-Norden manifold if 
\begin{equation*}
\underset{\xi _{1},\xi _{2},\xi _{3}}{\sigma }h\left( \left( \nabla _{\xi
_{1}}^{h}\tciLaplace \right) \xi _{2},\xi _{3}\right) =0,
\end{equation*}%
where $\sigma $ is the cyclic sum by three arguments.
\end{definition}

\begin{theorem}
\label{theolast}\cite{Salimov} Let $\left( \acute{N},\tciLaplace ,h\right) $
be non-integrable an almost anti-Hermitian manifold. Then the triple $\left( 
\acute{N},\tciLaplace ,h\right) $ is a quasi-K\"{a}hler-Norden if and only if%
\begin{equation*}
\left( \Phi _{\tciLaplace }h\right) \left( \xi _{1},\xi _{2},\xi _{3}\right)
+\left( \Phi _{\tciLaplace }h\right) \left( \xi _{2},\xi _{3},\xi
_{1}\right) +\left( \Phi _{\tciLaplace }h\right) \left( \xi _{3},\xi
_{1},\xi _{2}\right) =0.
\end{equation*}
\end{theorem}

As is known, the almost complex structure $\tciLaplace $ on an anti-K\"{a}%
hler manifold $\left( \acute{N},h\right) $ is always integrable.

We will recall notions related to an anti-Hermitian metric.

The covariant derivative of the metrics $h$ and $\hslash $ are defined by%
\begin{equation*}
\left( \bigtriangledown _{\xi _{3}}h\right) \left( \xi _{1},\xi _{2}\right)
=\xi _{3}h\left( \xi _{1},\xi _{2}\right) -h\left( \bigtriangledown _{\xi
_{3}}\xi _{1},\xi _{2}\right) -h\left( \xi _{1},\bigtriangledown _{\xi
_{3}}\xi _{2}\right)
\end{equation*}%
and%
\begin{equation*}
\left( \bigtriangledown _{\xi _{3}}\hslash \right) \left( \xi _{1},\xi
_{2}\right) =\xi _{3}\hslash \left( \xi _{1},\xi _{2}\right) -\hslash \left(
\bigtriangledown _{\xi _{3}}\xi _{1},\xi _{2}\right) -\hslash \left( \xi
_{1},\bigtriangledown _{\xi _{3}}\xi _{2}\right) .
\end{equation*}%
Clearly $\left( \bigtriangledown _{\xi _{3}}h\right) \left( \xi _{1},\xi
_{2}\right) =\left( \bigtriangledown _{\xi _{3}}h\right) \left( \xi _{2},\xi
_{1}\right) $ and $\left( \bigtriangledown _{\xi _{3}}\hslash \right) \left(
\xi _{1},\xi _{2}\right) =\left( \bigtriangledown _{\xi _{3}}\hslash \right)
\left( \xi _{2},\xi _{1}\right) $ due to symmetry of $h$ and $\hslash $. For
any linear connection $\bigtriangledown $, its $h$-conjugate connection $%
\bigtriangledown ^{\sharp }$ and its $\hslash $-conjugate connection $%
\bigtriangledown ^{\ddag }$ are defined by%
\begin{equation*}
\xi _{3}h\left( \xi _{1},\xi _{2}\right) =h\left( \bigtriangledown _{\xi
_{3}}\xi _{1},\xi _{2}\right) +h\left( \xi _{1},\bigtriangledown _{\xi
_{3}}^{\sharp }\xi _{2}\right)
\end{equation*}%
and%
\begin{equation*}
\xi _{3}\hslash \left( \xi _{1},\xi _{2}\right) =\hslash \left(
\bigtriangledown _{\xi _{3}}\xi _{1},\xi _{2}\right) +\hslash \left( \xi
_{1},\bigtriangledown _{\xi _{3}}^{\ddag }\xi _{2}\right) ,
\end{equation*}%
respectively. It can be easily checked that the following conditions are
satisfied 
\begin{equation*}
\left( \bigtriangledown ^{\sharp }\right) ^{\sharp }=\left( \bigtriangledown
^{\ddag }\right) ^{\ddag }=\left( \bigtriangledown ^{\tciLaplace }\right)
^{\tciLaplace }=\bigtriangledown ,
\end{equation*}%
\begin{equation*}
\bigtriangledown ^{\sharp }=\left( \bigtriangledown ^{\dag }\right)
^{\tciLaplace }=\left( \bigtriangledown ^{\tciLaplace }\right) ^{\dag },
\end{equation*}%
\begin{equation*}
\bigtriangledown ^{\ddag }=\left( \bigtriangledown ^{\sharp }\right)
^{\tciLaplace }=\left( \bigtriangledown ^{\tciLaplace }\right) ^{\sharp },
\end{equation*}%
\begin{equation*}
\bigtriangledown ^{\tciLaplace }=\left( \bigtriangledown ^{\sharp }\right)
^{\ddag }=\left( \bigtriangledown ^{\ddag }\right) ^{\sharp },
\end{equation*}%
which give that $\left( id,\sharp ,\ddag ,\tciLaplace \right) $ is a $4$%
-element Klein group action on the space of linear connections (also see 
\cite{Gezer}). Next, we will give some results without proof. These results
can be proven by following the proofs of Propositions \ref{pro3}, \ref{pro4}
and \ref{pro5}. Their proofs are used the anti-Hermitian metric $h$ and the
twin anti-Hermitian metric $\hslash $ instead of the Hermitian metric $g$
and the fundamental $2$-form $\omega $ and purity conditions.

\begin{proposition}
Let $\left( \tciLaplace ,h\right) $ be an almost anti-Hermitian manifold and
let $\bigtriangledown $ be a linear connection with torsion tensor $%
T^{\nabla }$ on $\acute{N}$. Let $\hslash $ be the twin anti-Hermitian
metric. Then, there exist the following expressions

$\left( i\right) $ Assume that $\left( \bigtriangledown ,\tciLaplace \right) 
$ is Codazzi-coupled. $\left( \bigtriangledown ^{\sharp },\hslash \right) $
is a quasi-statistical structure $\Leftrightarrow \left( \bigtriangledown
^{\sharp },h\right) $ is a quasi-statistical structure.

$\left( ii\right) $ Assume that $\left( \bigtriangledown ,\tciLaplace
\right) $ is Codazzi-coupled. $\left( \bigtriangledown ^{\ddag },\hslash
\right) $ is a quasi-statistical structure $\Leftrightarrow \left(
\bigtriangledown ^{\ddag },h\right) $ is a quasi-statistical structure.

$\left( iii\right) $ Assume that $\left( \bigtriangledown ^{\sharp
},\tciLaplace \right) $ is Codazzi-coupled. $\left( \bigtriangledown
,\hslash \right) $ is a quasi-statistical structure $\Leftrightarrow \left(
\bigtriangledown ,h\right) $ is a quasi-statistical structure.

$\left( iv\right) $ Assume that $\left( \bigtriangledown ^{\ddag
},\tciLaplace \right) $ is Codazzi-coupled. $\left( \bigtriangledown
,\hslash \right) $ is a quasi-statistical structure $\Leftrightarrow \left(
\bigtriangledown ,h\right) $ is a quasi-statistical structure.

$\left( v\right) $ Assume that $\left( \bigtriangledown ^{\ddag
},\tciLaplace \right) $ is Codazzi-coupled. $\left( \bigtriangledown
^{\tciLaplace },\hslash \right) $ is a quasi statistical structure $%
\Leftrightarrow \left( \bigtriangledown ,\hslash \right) $ is a quasi
statistical structure.

$\left( vi\right) $ Assume that $\left( \bigtriangledown ^{\sharp
},\tciLaplace \right) $ is Codazzi-coupled. $\left( \bigtriangledown
^{\tciLaplace },h\right) $ is a quasi-statistical structure $\Leftrightarrow
\left( \bigtriangledown ,h\right) $ is a quasi-statistical structure.
\end{proposition}

\begin{proposition}
\label{pro12}Let $\left( \acute{N},\tciLaplace ,h\right) $ be an
anti-Hermitian manifold and let $\bigtriangledown $ be a linear connection
with torsion tensor $T^{\nabla }$ on $\acute{N}$. Let $\hslash $ be the twin
anti-Hermitian metric. Then, the following expressions hold

$\left( i\right) $ $\left( \bigtriangledown ^{\tciLaplace },\hslash \right) $
is a quasi-statistical structure if and only if $\left( \bigtriangledown
,h\right) $ is a quasi-statistical structure.

$\left( ii\right) $ $\left( \bigtriangledown ,\hslash \right) $ is a
quasi-statistical structure if and only if $\left( \bigtriangledown
^{\tciLaplace },h\right) $ is a quasi-statistical structure.

$\left( iii\right) $ $\left( \bigtriangledown ^{\ddag },\hslash \right) $ is
a quasi-statistical structure if and only if $\left( \bigtriangledown
^{\sharp },h\right) $ is a quasi-statistical structure.

$\left( iv\right) $ $\left( \bigtriangledown ^{\sharp },\hslash \right) $ is
a quasi-statistical structure if and only if $\left( \bigtriangledown
^{\ddag },h\right) $ is a quasi-statistical structure.
\end{proposition}

\begin{corollary}
\label{cor7}Let $\acute{N}$ be a manifold equipped with an anti-Hermitian
metric $h$, a linear connection $\bigtriangledown $ with torsion tensor $%
T^{\bigtriangledown }$ and the twin anti-Hermitian metric $\hslash $. Denote
by $\bigtriangledown ^{\sharp },\bigtriangledown ^{\ddag }$and $%
\bigtriangledown ^{\tciLaplace }$, respectively, $h$-conjugation, $\hslash $%
-conjugation and $\tciLaplace $-conjugate transformations of an arbitrary
linear connection $\bigtriangledown .$ From $4$-element Klein group action
on the space of linear connections, we have

$\left( i\right) $ $\left( \bigtriangledown ,h\right) $ is a
quasi-statistical structure if and only if the linear connection $%
\bigtriangledown ^{\sharp }$ is torsion-free.

$\left( ii\right) $ $\left( \bigtriangledown ^{\sharp },h\right) $ is a
quasi-statistical structure if and only if the linear connection $%
\bigtriangledown $ is torsion-free.

$\left( iii\right) $ $\left( \bigtriangledown ^{\tciLaplace },h\right) $ is
a quasi-statistical structure if and only if the linear connection $%
\bigtriangledown ^{\ddag }$ is torsion-free.

$\left( iv\right) $ $\left( \bigtriangledown ^{\ddag },h\right) $ is a
quasi-statistical structure if and only if the linear connection $%
\bigtriangledown ^{\tciLaplace }$ is torsion-free.
\end{corollary}

From Proposition \ref{pro12} and Corollary \ref{cor7}, we have the following
result.

\begin{proposition}
Let $\left( \acute{N},\tciLaplace ,h\right) $ be an almost anti-Hermitian
manifold, $\bigtriangledown $ be an arbitrary linear connection, $%
\bigtriangledown ^{\sharp }$ be the $h$-conjugate connection of $%
\bigtriangledown $ and $\bigtriangledown ^{\ddag }$be the $\hslash $%
-conjugate connection of $\bigtriangledown $. Then, there exist the below
expressions

$\left( i\right) $ $d^{\bigtriangledown }\hslash =0\Leftrightarrow
T^{\bigtriangledown ^{\ddag }}=0\Leftrightarrow d^{\bigtriangledown ^{\sharp
}}\tciLaplace =0\Leftrightarrow d^{\bigtriangledown ^{\tciLaplace }}h=0$;

$\left( ii\right) $ $d^{\bigtriangledown ^{\sharp }}\hslash
=0\Leftrightarrow T^{\bigtriangledown ^{\tciLaplace }}=0\Leftrightarrow
d^{\bigtriangledown }\tciLaplace =0\Leftrightarrow d^{\bigtriangledown
^{\ddag }}h=0$;

$\left( iii\right) $ $d^{\bigtriangledown ^{\ddag }}\hslash
=0\Leftrightarrow T^{\bigtriangledown }=0\Leftrightarrow d^{\bigtriangledown
^{\tciLaplace }}\tciLaplace =0\Leftrightarrow d^{\bigtriangledown ^{\sharp
}}h=0$;

$\left( iv\right) $ $d^{\bigtriangledown ^{\tciLaplace }}\hslash
=0\Leftrightarrow T^{\bigtriangledown ^{\sharp }}=0\Leftrightarrow
d^{\bigtriangledown ^{\ddag }}\tciLaplace =0\Leftrightarrow
d^{\bigtriangledown }h=0$.
\end{proposition}

\begin{proposition}
\label{pro14}Let $\left( \acute{N},\tciLaplace ,h\right) $ be an almost
anti-Hermitian manifold and $\bigtriangledown $ be a linear connection with
torsion tensor $T^{\nabla }$ on $\acute{N}$. If $\left( \bigtriangledown
,h\right) $ is a quasi-statistical structure, that is, $d^{\bigtriangledown
}h=0$, we get%
\begin{eqnarray*}
\left( \Phi _{\tciLaplace }h\right) \left( \xi _{1},\xi _{2},\xi _{3}\right)
&=&\left( \bigtriangledown _{\xi _{2}}h\right) \left( \tciLaplace \xi
_{1},\xi _{3}\right) -\left( \bigtriangledown _{\xi _{2}}h\right) \left( \xi
_{1},\tciLaplace \xi _{3}\right) +h\left( \left( \bigtriangledown _{\xi
_{2}}\tciLaplace \right) \xi _{1},\xi _{3}\right) \\
&&+h\left( \xi _{2},\left( \bigtriangledown _{\xi _{3}}\tciLaplace \right)
\xi _{1}\right) -h\left( \xi _{2},\left( \bigtriangledown _{\xi
_{1}}\tciLaplace \right) \xi _{3}\right) \\
&&+h\left( \xi _{2},T^{\bigtriangledown }\left( \tciLaplace \xi _{1},\xi
_{3}\right) -\tciLaplace T^{\bigtriangledown }\left( \xi _{1},\xi
_{3}\right) \right) .
\end{eqnarray*}
\end{proposition}

\begin{proof}
Using the definition of the Tachibana operator $\Phi _{\tciLaplace }$, we
have%
\begin{eqnarray*}
\left( \Phi _{\tciLaplace }h\right) \left( \xi _{1},\xi _{2},\xi _{3}\right)
&=&\tciLaplace \xi _{1}h\left( \xi _{2},\xi _{3}\right) -\xi _{1}h\left(
\tciLaplace \xi _{2},\xi _{3}\right) \\
&&+h\left( \left( L_{\xi _{2}}\tciLaplace \right) \xi _{1},\xi _{3}\right)
+h\left( \xi _{2},\left( L_{\xi _{3}}\tciLaplace \right) \xi _{1}\right) ,
\end{eqnarray*}%
where $\left( L_{\xi _{1}}\tciLaplace \right) \xi _{2}=\left[ \xi
_{1},\tciLaplace \xi _{2}\right] -\tciLaplace \left[ \xi _{1},\xi _{2}\right]
$. Then, we obtain%
\begin{eqnarray*}
&&\left( \Phi _{\tciLaplace }h\right) \left( \xi _{1},\xi _{2},\xi
_{3}\right) \\
&=&\left( \bigtriangledown _{\tciLaplace \xi _{1}}h\right) \left( \xi
_{2},\xi _{3}\right) -\left( \bigtriangledown _{\xi _{1}}h\right) \left( \xi
_{2},\tciLaplace \xi _{3}\right) +h\left( T^{\bigtriangledown }\left(
\tciLaplace \xi _{1},\xi _{2}\right) ,\xi _{3}\right) \\
&&-h\left( T^{\bigtriangledown }\left( \xi _{1},\xi _{2}\right) ,\tciLaplace
\xi _{3}\right) +h\left( \left( \bigtriangledown _{\xi _{2}}\tciLaplace
\right) \xi _{1},\xi _{3}\right) +h\left( \xi _{2},\left( \bigtriangledown
_{\xi _{3}}\tciLaplace \right) \xi _{1}\right) \\
&&-h\left( \xi _{2},\left( \bigtriangledown _{\xi _{1}}\tciLaplace \right)
\xi _{3}\right) +h\left( \xi _{2},T^{\bigtriangledown }\left( \tciLaplace
\xi _{1},\xi _{3}\right) -\tciLaplace T^{\bigtriangledown }\left( \xi
_{1},\xi _{3}\right) \right) .
\end{eqnarray*}%
Since $\left( \bigtriangledown ,h\right) $ is a quasi-statistical structure,
it follows that%
\begin{eqnarray*}
\left( \Phi _{\tciLaplace }h\right) \left( \xi _{1},\xi _{2},\xi _{3}\right)
&=&\left( \bigtriangledown _{\xi _{2}}h\right) \left( \tciLaplace \xi
_{1},\xi _{3}\right) -\left( \bigtriangledown _{\xi _{2}}h\right) \left( \xi
_{1},\tciLaplace \xi _{3}\right) \\
&&+h\left( \left( \bigtriangledown _{\xi _{2}}\tciLaplace \right) \xi
_{1},\xi _{3}\right) +h\left( \xi _{2},\left( \bigtriangledown _{\xi
_{3}}\tciLaplace \right) \xi _{1}\right) \\
&&-h\left( \xi _{2},\left( \bigtriangledown _{\xi _{1}}\tciLaplace \right)
\xi _{3}\right) +h\left( \xi _{2},T^{\bigtriangledown }\left( \tciLaplace
\xi _{1},\xi _{3}\right) -\tciLaplace T^{\bigtriangledown }\left( \xi
_{1},\xi _{3}\right) \right) .
\end{eqnarray*}
\end{proof}

Let us sign that $h\left( T^{\nabla }\left( \tciLaplace \xi _{1},\xi
_{2}\right) ,\xi _{3}\right) =T^{\nabla }\left( \tciLaplace \xi _{1},\xi
_{2},\xi _{3}\right) $ and $h\left( \left( \nabla _{\xi _{1}}\tciLaplace
\right) \xi _{2},\xi _{3}\right) =B\left( \xi _{1},\xi _{2},\xi _{3}\right) $%
. Hence, we say that if $T^{\nabla }\left( \tciLaplace \xi _{1},\xi _{3},\xi
_{2}\right) =-B\left( \xi _{2},\xi _{3},\xi _{1}\right) $, then we have $%
T^{\nabla }\left( \tciLaplace \xi _{1},\xi _{2}\right) =-T^{\nabla }\left(
\xi _{1},\tciLaplace \xi _{2}\right) $. Hence, we are ready to give the
second main theorem of this paper.

\begin{theorem}
\label{teo5}Let $\left( \acute{N},\tciLaplace ,h\right) $ be an almost
anti-Hermitian manifold and $\bigtriangledown $ be a linear connection with
torsion tensor $T^{\nabla }$ on $\acute{N}$. Suppose that $%
d^{\bigtriangledown }\tciLaplace =0$ and $d^{\bigtriangledown }h=0$. Then,
the triple $\left( \acute{N},\tciLaplace ,g\right) $ is an anti-K\"{a}hler
manifold if and only if the condition $T^{\nabla }\left( \tciLaplace \xi
_{1},\xi _{3},\xi _{2}\right) =-B\left( \xi _{2},\xi _{3},\xi _{1}\right) $
for any vector fields $\xi _{1}$, $\xi _{2}$, $\xi _{3}$ on $\acute{N}$
holds.
\end{theorem}

\begin{proof}
From the Proposition \ref{pro14}, if $\left( \nabla ,h\right) $ is a
quasi-statistical structure, then we get 
\begin{eqnarray*}
\left( \Phi _{\tciLaplace }h\right) \left( \xi _{1},\xi _{2},\xi _{3}\right)
&=&\left( \bigtriangledown _{\xi _{2}}h\right) \left( \tciLaplace \xi
_{1},\xi _{3}\right) -\left( \bigtriangledown _{\xi _{2}}h\right) \left( \xi
_{1},\tciLaplace \xi _{3}\right) +h\left( \left( \bigtriangledown _{\xi
_{2}}\tciLaplace \right) \xi _{1},\xi _{3}\right) \\
&&+h\left( \xi _{2},\left( \bigtriangledown _{\xi _{3}}\tciLaplace \right)
\xi _{1}\right) -h\left( \xi _{2},\left( \bigtriangledown _{\xi
_{1}}\tciLaplace \right) \xi _{3}\right) \\
&&+h\left( \xi _{2},T^{\bigtriangledown }\left( \tciLaplace \xi _{1},\xi
_{3}\right) -\tciLaplace T^{\bigtriangledown }\left( \xi _{1},\xi
_{3}\right) \right) .
\end{eqnarray*}%
Considering the condition $d^{\bigtriangledown }\tciLaplace =0$, we have%
\begin{equation*}
\left( \Phi _{\tciLaplace }h\right) \left( \xi _{1},\xi _{2},\xi _{3}\right)
=T^{\bigtriangledown }\left( \tciLaplace \xi _{1},\xi _{3},\xi _{2}\right)
+B\left( \xi _{2},\xi _{3},\xi _{1}\right) ,
\end{equation*}%
from which we immediately say that the triple $\left( \acute{N},\tciLaplace
,g\right) $ is an anti-K\"{a}hler manifold if and only if the condition $%
T^{\bigtriangledown }\left( \tciLaplace \xi _{1},\xi _{3},\xi _{2}\right)
=-B\left( \xi _{2},\xi _{3},\xi _{1}\right) $ holds.
\end{proof}

\begin{remark}
The Theorem \ref{teo5} says that for any quasi-statistical manifold $\left( 
\acute{N},h,\nabla \right) $, if the almost complex structure $\tciLaplace $
is $d^{\nabla }$-closed and the condition $T^{\bigtriangledown }\left(
\tciLaplace \xi _{1},\xi _{3},\xi _{2}\right) +B\left( \xi _{2},\xi _{3},\xi
_{1}\right) =0$ is satisfied, $\left( \acute{N},h,\nabla ,\tciLaplace
\right) $ is anti-K\"{a}hler manifold. By taking any linear connection $%
\nabla $ with torsion tensor $T^{\nabla }$ instead of Levi-Civita connection 
$\nabla ^{h}$ of $h$ or torsion-free linear connection, it is also possible
to make a classification for anti-K\"{a}hler manifolds.
\end{remark}

Let $\left( \acute{N},\tciLaplace ,h\right) $ be an almost anti-Hermitian
manifold and $\bigtriangledown $ be a linear connection with torsion tensor $%
T^{\nabla }$ on $\acute{N}$. If $\left( \bigtriangledown ,h\right) $ is a
quasi-statistical structure, with help of Proposition \ref{pro14} we obtain%
\begin{eqnarray*}
&&\left( \Phi _{\tciLaplace }h\right) \left( \xi _{1},\xi _{2},\xi
_{3}\right) +\left( \Phi _{\tciLaplace }h\right) \left( \xi _{2},\xi
_{3},\xi _{1}\right) +\left( \Phi _{\tciLaplace }h\right) \left( \xi
_{3},\xi _{1},\xi _{2}\right) \\
&=&h\left( \xi _{2},T^{\bigtriangledown }\left( \tciLaplace \xi _{1},\xi
_{3}\right) \right) +h\left( \xi _{3},T^{\bigtriangledown }\left(
\tciLaplace \xi _{2},\xi _{1}\right) \right) +h\left( \xi
_{1},T^{\bigtriangledown }\left( \tciLaplace \xi _{3},\xi _{2}\right) \right)
\\
&&+h\left( \xi _{2},\left( \bigtriangledown _{\xi _{3}}\tciLaplace \right)
\xi _{1}\right) +h\left( \xi _{3},\left( \bigtriangledown _{\xi
_{1}}\tciLaplace \right) \xi _{2}\right) +h\left( \xi _{1},\left(
\bigtriangledown _{\xi _{2}}\tciLaplace \right) \xi _{3}\right) .
\end{eqnarray*}%
Hence, the last equality and Theorem \ref{theolast} give the following
result.

\begin{theorem}
\label{cor8}Let $\left( \acute{N},\tciLaplace ,h\right) $ be an almost
anti-Hermitian manifold and $\bigtriangledown $ be a linear connection with
torsion tensor $T^{\nabla }$ on $\acute{N}$. Under the assumption that $%
\left( \bigtriangledown ,h\right) $ is a quasi-statistical structure, the
triple $\left( \acute{N},\tciLaplace ,h\right) $ is a quasi-K\"{a}%
hler-Norden manifold if and only if%
\begin{eqnarray*}
&&h\left( \xi _{2},\left( \bigtriangledown _{\xi _{3}}\tciLaplace \right)
\xi _{1}\right) +h\left( \xi _{3},\left( \bigtriangledown _{\xi
_{1}}\tciLaplace \right) \xi _{2}\right) +h\left( \xi _{1},\left(
\bigtriangledown _{\xi _{2}}\tciLaplace \right) \xi _{3}\right) \\
&=&-\left( h\left( \xi _{2},T^{\bigtriangledown }\left( \tciLaplace \xi
_{1},\xi _{3}\right) \right) +h\left( \xi _{3},T^{\bigtriangledown }\left(
\tciLaplace \xi _{2},\xi _{1}\right) \right) +h\left( \xi
_{1},T^{\bigtriangledown }\left( \tciLaplace \xi _{3},\xi _{2}\right)
\right) \right) .
\end{eqnarray*}
\end{theorem}

\section{Conclusion}

It is known that the triple $\left( \acute{N},g,\tciLaplace \right) $ is a K%
\"{a}hler manifold with an almost complex structure $\tciLaplace $ and a
pseudo-Riemannian metric $g$ if and only if the almost complex structure $%
\tciLaplace $ is parallel under the Levi-Civita connection $\nabla ^{g}$ of $%
g$. In \cite{Fei}, using the Codazzi couplings of $\nabla $ with a
pseudo-Riemannian metric $g$ and an almost complex structure $\tciLaplace $,
the authors give a new alternative classification for K\"{a}hler manifold,
where $\nabla $ is any linear connection with torsion-free. In this paper,
for any linear connection $\nabla $ with torsion tensor $T^{\nabla }$, we
have proven that under the assumption that $d^{\nabla }g=0$, $d^{\nabla
}\tciLaplace =0$ and $T^{\nabla }\left( \tciLaplace \xi _{1},\xi _{2}\right)
=-T^{\nabla }\left( \xi _{1},\tciLaplace \xi _{2}\right) $, the almost
complex structure $\tciLaplace $ is integrable and $\omega $ is closed.
Hence, the almost Hermitian manifold $\left( \acute{N},g,\tciLaplace \right) 
$ rises to a K\"{a}hler manifold. This shows us that it is not necessary for
the connection to be torsion-free to make such a classification. Moreover,
this paper shows that under the certain conditions the anti-K\"{a}hler
manifolds can be characterized by taking any linear connection $\nabla $
with torsion tensor $T^{\nabla }$ instead of the Levi-Civita connection $%
\nabla ^{h}$ of a pseudo-Riemannian metric $h$ or torsion-free linear
connection. Consequently, the paper gives new classifications for the K\"{a}%
hler, anti-K\"{a}hler and quasi-K\"{a}hler-Norden manifolds.

\end{document}